\documentclass[oneside,english,reqno]{amsart}
\usepackage{ae,aecompl}
\usepackage[T1]{fontenc}
\usepackage[latin9]{inputenc}
\usepackage{geometry}
\usepackage{mathtools}
\usepackage{amstext}
\usepackage{amsthm}
\usepackage{amssymb}
\usepackage{graphicx}
\PassOptionsToPackage{normalem}{ulem}
\usepackage{ulem}

\makeatletter
\theoremstyle{plain}
\newtheorem{thm}{\protect\theoremname}[section]
\theoremstyle{plain}
\newtheorem{fact}[thm]{\protect\factname}
\theoremstyle{plain}
\newtheorem{question}[thm]{\protect\questionname}
\theoremstyle{plain}
\newtheorem{lem}[thm]{\protect\lemmaname}
\theoremstyle{plain}
\newtheorem{conjecture}[thm]{\protect\conjecturename}
\theoremstyle{definition}
\newtheorem{defn}[thm]{\protect\definitionname}
\theoremstyle{definition}
\newtheorem{example}[thm]{\protect\examplename}
\theoremstyle{plain}
\newtheorem{prop}[thm]{\protect\propositionname}
\theoremstyle{plain}
\newtheorem{cor}[thm]{\protect\corollaryname}
\theoremstyle{definition}
\newtheorem{condition}[thm]{\protect\conditionname}
\theoremstyle{definition}
\newtheorem{problem}[thm]{\protect\problemname}

\usepackage{enumitem} % Enumerate configuration
\usepackage[colorlinks = true,
            linkcolor = blue,
            urlcolor  = blue,
            citecolor = blue,
            anchorcolor = blue]{hyperref} % Hyper Links
\usepackage{fancyhdr} % Header/Footer

\allowdisplaybreaks

\emergencystretch=\maxdimen
\newcommand{\zGetCorrectThmNumber}{\arabic{enumi}.\arabic{enumii}.\arabic{thm}}

\hypersetup{pdfborder={0 0 0}}

\newcounter{zLinkCnt}
\newcommand{\zlabel}[2]{
\setcounter{thm}{\value{thm}+1}
{
\renewcommand\thezLinkCnt{#1 \zGetCorrectThmNumber}
\refstepcounter{zLinkCnt}
\label{z#2}
}
\setcounter{thm}{\value{thm}-1}
}

\newcommand{\zgoto}[1]{\ref{z#1}}

\newcommand{\zTableOfContents}{
\renewcommand\contentsname{Table Of Contents}
\tableofcontents
}

\setenumerate[enumerate,1]
{labelindent=\parindent,itemindent=1cm,leftmargin=8pt,labelsep=*, label*=\arabic*,font=\huge\textbf}

\setenumerate[enumerate,2]
{labelindent=3pt,itemindent=2pc,leftmargin=0pt,labelsep=*, label*=.\arabic*,font=\large\textbf}

\setenumerate[enumerate,3]{label=\arabic*.}

\newcommand{\zChapterHeader}[1]{{\huge \textbf {#1}} \bigskip}

\newcommand{\zChapter}[1]{
\item \zChapterHeader{#1}
\addcontentsline{toc}{chapter}{\textbf{\arabic{enumi} #1}}}

\newcommand{\zEndChapter}{
\ifnum \value{enumii}>0 \end{enumerate} \setcounter{enumii}{0} \fi
\newpage
}

\newcommand{\zSection}[1]{
\ifnum \value{enumii}=0 \begin{enumerate} \fi
\medskip\bigskip \item {\large \textbf {#1}} \bigskip \refstepcounter{section}
\addcontentsline{toc}{section}{\arabic{enumi}.\arabic{enumii} #1}}

\newcommand{\zBibRef}[1]{[\hyperlink{zBibRef#1}{#1}]}

\newcommand{\zBibSrc}[2]{\item [{[#1]}] {\hypertarget{zBibRef#1} #2}}

\makeatother

\usepackage{babel}
\providecommand{\conditionname}{Condition}
\providecommand{\conjecturename}{Conjecture}
\providecommand{\corollaryname}{Corollary}
\providecommand{\definitionname}{Definition}
\providecommand{\examplename}{Example}
\providecommand{\factname}{Fact}
\providecommand{\lemmaname}{Lemma}
\providecommand{\problemname}{Problem}
\providecommand{\propositionname}{Proposition}
\providecommand{\questionname}{Question}
\providecommand{\theoremname}{Theorem}

\begin{document}
\includegraphics[scale=0.8]{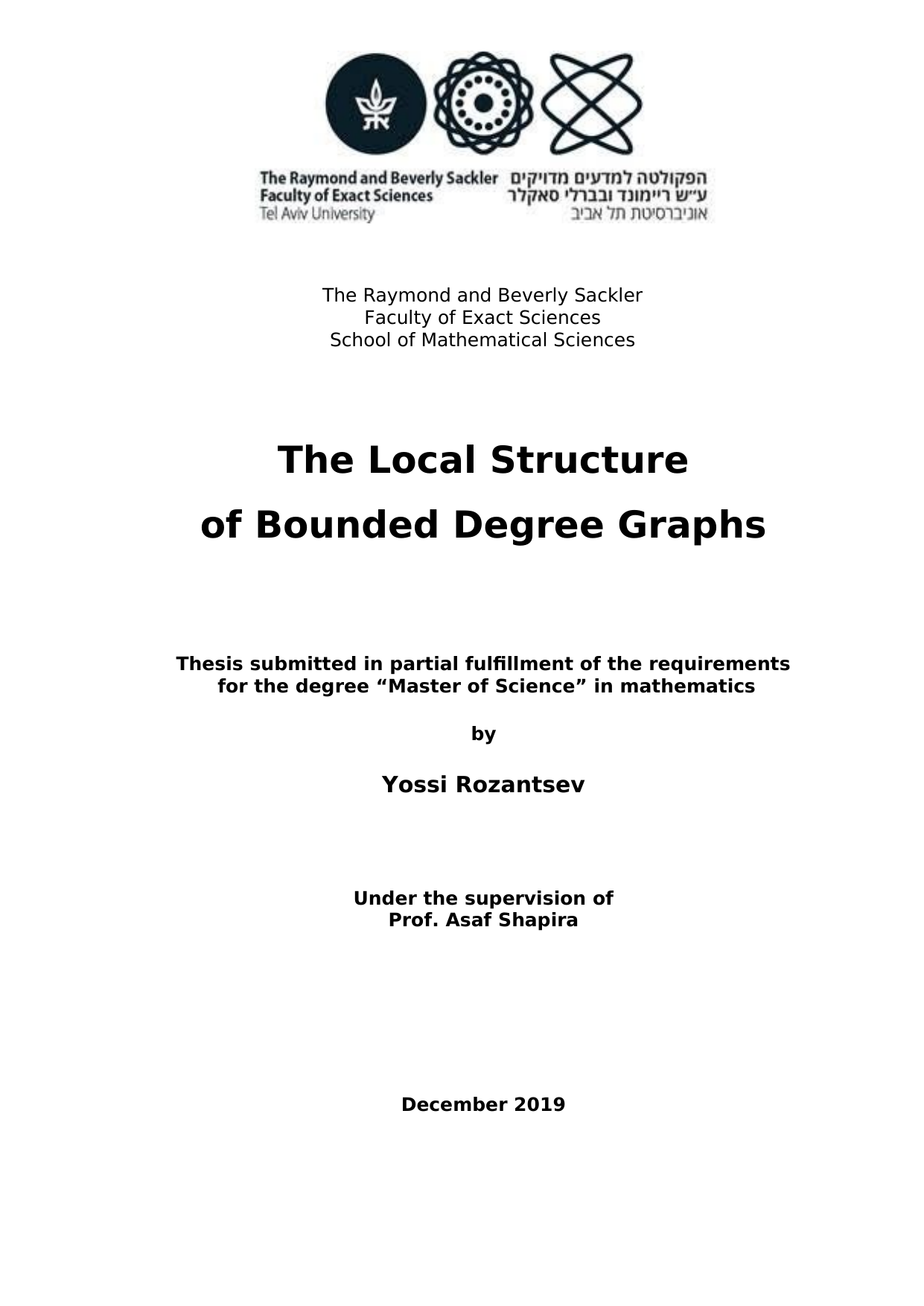}

\noindent \noindent \zChapterHeader{Acknowledgements}\\

\noindent First and foremost, I would like to thank my adviser, Asaf
Shapira, for accepting me as his student. Without any doubt this has
been one of the most special experiences of my life, and I'm glad
that it was done under his unique guidance. His endless patience,
and ability to always see more and beyond every small detail meant
a lot to me, and it was my honor to work with such a great mathematician.
\\
\\
I would like to thank my friends and companions: Ziv Greenhut, Gal
Ordo, Gur Lifshitz, Saar Diskin and Emanuel Segal for being there
for me during hard times and for being with me on this journey.\\
\\
I would like to thank my family for giving me the foundation which
brought me to this point, and for believing and supporting me throughout
these last few years.\\
\\
Finally, I would like to thank LyX/TeX for teaching me that there
is always hope, even when it seems that it's impossible to add an
empty line after a paragraph without breaking the entire PDF or crashing
the operating system.

\noindent \begin{flushleft}
\zEndChapter
\par\end{flushleft}

\noindent \noindent \begin{flushleft}
\zChapterHeader{Abstract}
\par\end{flushleft}

\noindent Let $G=(V,E)$ be a simple graph with maximum degree $d$.
For an integer $k\in\mathbb{N}$, the $k$-disc of a vertex $v\in V$
is defined as the rooted subgraph of $G$ that is induced by all vertices
whose distance to $v$ is at most $k$. The $k$-disc frequency distribution
vector of $G$, denoted by $\text{freq}_{k}(G)$, is a vector indexed
by all isomorphism types of rooted $k$-discs. For each such isomorphism
type $\Gamma$, the corresponding entry in $\text{freq}_{k}(G)$ counts
the fraction of vertices in $V$ that have a $k$-disc isomorphic
to $\Gamma$. In a sense, $\text{freq}_{k}(G)$ is one way to represent
the ``local structure'' of $G$.\\
\\
The graph $G$ can be arbitrarily large, and so a natural question
is whether given $\text{freq}_{k}(G)$ it is possible to construct
a small graph $H$, whose size is independent of $|V|$, such that
$H$ has a similar local structure. N. Alon proved that for any $\epsilon>0$
there always exists a graph $H$ whose size is independent of $|V|$
and whose frequency vector satisfies $||\text{freq}_{k}(G)-\text{freq}_{k}(H)||_{1}\le\epsilon$.
However, his proof is only existential and does not imply that there
is a deterministic algorithm to construct such a graph $H$. He gave
the open problem of finding an explicit deterministic algorithm that
finds $H$, or proving that no such algorithm exists.\\
\\
A possible approach to showing that there is no deterministic algorithm
that solves the problem is by reduction from a different undecidable
problem. This approach was used by P. Winkler to prove a similar theorem
- given a set $\varPhi$ of $k$-discs of a directed edge-colored
graph it is not possible to determine whether there exists a graph
whose set of $k$-discs is exactly $\varPhi$. The reduction was done
from a variant of the Post Correspondence Problem (PCP), which is
known to be undecidable. It is therefore interesting to examine the
directed edge-colored variant of Alon's question, and it's connection
to PCP.\\

\noindent Our main result is that Alon's problem is undecidable if
and only if the much more general problem (involving directed edges
and edge colors) is undecidable. We also prove that both problems
are decidable for the special case when $G$ is a path. We show that
the local structure of any directed edge-colored path $G$ can be
approximated by a suitable fixed-size directed edge-colored path $H$
and we give explicit bound on the size of $H$.

\noindent \begin{flushleft}
\zEndChapter
\zTableOfContents
\zEndChapter
% Chapter enumeration begins in the first indexed chapter.
\par\end{flushleft}

\noindent \noindent \begin{enumerate} % Begin the chapter enumeration
\zChapter{Introduction}

\noindent Let $d\ge2$ and $k\ge1$ be fixed integers. A simple graph
$G$ is a finite, unweighted, undirected graph containing no loops
or multiple edges. We write $G=(V,E)$, where $V=V(G)$ is a finite
set of vertices, and $E=E(G)$ is the set of edges. Throughout the
thesis, we will assume $G$ to be a $d$-bounded degree graph, that
is, the maximum degree of a vertex in $G$ is upper bounded by $d$.
Given two different vertices $u,v\in V$ let $\text{dist}_{G}(u,v)$
be the length of the shortest path between $u$ and $v$.

\noindent For any $v\in V$, the \textbf{$k$-disc} of $v$, denoted
by $\text{disc}_{k}(v)$ or $\text{disc}_{k}(G,v)$ is defined as
the rooted subgraph in $G$ that is induced by the vertices that are
at distance at most $k$ to $v$ in $G$. Two $k$-discs are isomorphic
if and only if there exists a root-preserving graph isomorphism between
them (a graph isomorphism that identifies the roots). We denote the
set of all non-isomorphic $d$-bounded degree rooted graphs with radius
at most $k$ by $\mathcal{L}(d,k)$.

\noindent \zlabel{Fact}{DiscSize}
\begin{fact}
\noindent Let $v\in V$ be a vertex, then $|\text{disc}_{k}(v)|\le2d^{k}$.
In particular $\mathcal{L}(d,k)$ is finite.
\end{fact}

\noindent 
\noindent \textbf{Proof} Since $d$ is finite, the amount of vertices
in $\text{disc}_{k}(v)$ is at most $1+d+...+d^{k}\le2d^{k}$. There
is only a finite amount of simple $d$-bounded degree graphs on $2d^{k}$
vertices, and so $\mathcal{L}(d,k)$ is finite.\\
\\
We denote the size of $\mathcal{L}(d,k)$ by $L\coloneqq L(d,k)$
and write $\mathcal{L}(d,k)=\{\Gamma_{1},...,\Gamma_{L}\}$.

\noindent The $k$-disc \textbf{count vector} $\text{cnt}_{k}(G)$
of a graph $G$ is an $L$-dimensional vector where the $i$-th entry
counts the number of $k$-discs in $G$ that are isomorphic to $\Gamma_{i}\in\mathcal{L}(d,k)$.
Given a $k$-disc isomorphism type $\Gamma$, $\text{cnt}_{k}(G,\Gamma)$
is defined as the entry in $\text{cnt}_{k}(G)$ that corresponds to
$\Gamma$.

\noindent The $k$-disc \textbf{frequency distribution vector} \textbf{(FDV)}
of $G$, denoted by $\text{freq}_{k}(G)$, is the vector where the
$i$-th entry counts the fraction of $k$-discs in $G$ that are isomorphic
to $\Gamma_{i}\in\mathcal{L}(d,k)$, or equivalently $\text{freq}_{k}(G)\coloneqq\text{cnt}_{k}(G)/|V(G)|$.

\noindent Given a $k$-disc isomorphism type $\Gamma$, $\text{freq}_{k}(G,\Gamma)$
is defined as the entry in $\text{freq}_{k}(G)$ that corresponds
to $\Gamma$.\\

\noindent The main question, as given by N. Alon in \zBibRef{BER},
is whether it is possible to construct a small graph $H$ that has
approximately the same local structure as an arbitrarily large simple
graph $G$ whose degree is bounded by $d$.

\zlabel{Question}{Main}
\begin{question}
\noindent (N. Alon)

\noindent Let $d\ge2,k\ge1,\epsilon\in(0,1)$. Is there a computable
function $f\coloneqq f(d,k,\epsilon)$ such that for any simple $d$-bounded
graph $G$ there is a simple graph $H$ such that
\[
||\text{freq}_{k}(G)-\text{freq}_{k}(H)||_{1}\le\epsilon\quad\wedge\quad|V(H)|\le f(d,k,\epsilon)
\]
\end{question}

\noindent 
\noindent The motivation behind finding such an approximation is that
any algorithm that only uses the local structure of a graph will behave
similarly on $G$ and $H$. This is interesting in the context of
property testing in the bounded degree graph model, as introduced
by Goldreich and Ron \zBibRef{GR}, where we are given access to the
adjacency lists of vertices in a graph $G$ with maximum degree $d$,
and the goal is to distinguish between graphs with a given property
$\Pi$ and graphs that are $\epsilon$-far from having $\Pi$, that
is, graphs in which at least $\epsilon d|G|$ edges need to be changed
for the graph to have the property $\Pi$. There are many property
testers in this model that only depend on the local structure of graphs.
For example, all minor-closed properties can be tested this way (see
\zBibRef{BSS} and \zBibRef{HKN}).

\noindent If we look at dense graphs, instead of bounded degree graphs,
and replace $k$-discs with induced subgraphs of size $k$, then it
is possible to find a small graph $H$ whose local structure is close
to that of $G$. This follows from the regularity lemma \zBibRef{REG},
which provides a constant size weighted graph that captures the local
structure of $G$.

\noindent \zSection{Known Results}

\noindent It was sketched by N. Alon that there is a well defined
function $f$ (not necessarily computable) that satisfies the required
condition (see \zBibRef{AL} Proposition $19.10$). We give the full
proof in the following lemma.

\noindent \zlabel{Lemma}{FiniteCover}
\begin{lem}
\noindent Let $d\ge2,k\ge1,\epsilon>0$. Then, there is a \textbf{finite}
set of simple graphs $W$ such that 
\[
|W|\le\left(\frac{2L(d,k)}{\epsilon}\right)^{L(d,k)}
\]

\noindent And for any simple graph $G$, there is a graph $H\in W$
with $||\text{freq}_{k}(G)-\text{freq}_{k}(H)||_{1}\le\epsilon$.

\noindent In particular, $f=\max_{H\in W}|V(H)|$ satisfies the condition
of Alon's question.
\end{lem}

\noindent \textbf{Proof} We know by \zgoto{DiscSize} that $L(d,k)$
is finite. We denote $L(d,k)$ by $n$ and define the following set
$X\subseteq[0,1]^{n}:$

\noindent 
\[
X=\{1\cdot\frac{\epsilon}{2n},2\cdot\frac{\epsilon}{2n},...,\lfloor\frac{2n}{\epsilon}\rfloor\cdot\frac{\epsilon}{2n}\}^{n}
\]

\noindent The set $X$ approximates any vector $v=(v_{1},...,v_{n})\in[0,1]^{n}$
up to an error of $\frac{\epsilon}{2n}$ per coordinate, namely
\[
||v-X||_{1}\le n\cdot\frac{\epsilon}{2n}=\frac{\epsilon}{2}
\]

\noindent We start with an empty set $W=\emptyset$, and for each
$x\in X$, if there is a simple graph $H=H(x)$ with ${||\text{freq}_{k}(H)-x||_{1}\le\frac{\epsilon}{2}}$,
we add $H$ to $W$. The size of $W$ in this case is at most $|X|$.
Moreover, we have
\[
|W|\le|X|\le\left(\lfloor\frac{2n}{\epsilon}\rfloor\right)^{n}\le\left(\frac{2n}{\epsilon}\right)^{n}=\left(\frac{2L(d,k)}{\epsilon}\right)^{L(d,k)}
\]

\noindent Finally, let $G$ be a simple graph. The $k$-disc frequency
distribution vector of $G$ is a vector in $[0,1]^{n}$, and so there
is an $x\in X$ such that $||\text{freq}_{k}(G)-x||_{1}\le\frac{\epsilon}{2}$.
In particular $H=H(x)$ is a well defined graph that is part of $W$
(possibly even $H=G$) and then

\noindent 
\begin{align*}
||\text{freq}_{k}(G)-\text{freq}_{k}(H)||_{1} & \le||\text{freq}_{k}(G)-x||_{1}+||x-\text{freq}_{k}(H)||_{1}\le\frac{\epsilon}{2}+\frac{\epsilon}{2}=\epsilon
\end{align*}

\noindent With $W$ being a finite set which instills an ``$\epsilon$-approximation''
of the local structure of all simple graphs, we conclude that $f=\max_{H\in W}|V(H)|$
satisfies the conditions of the main question.\qed\\
\\
In other words, there is a finite set of simple graphs $W$ such that
the local structure of any graph is ``approximated'' by a graph
in $W$. We have an upper bound on the minimum size of $W$, but the
proof does not give any information regarding which specific graphs
are in this set, or how many vertices they have.\\
\\
Partial progress towards answering the main question was done by Fichtenberger,
Peng and Sohler \zBibRef{FPS}. They have shown that the computable
function
\[
f(d,k,\epsilon)=36\frac{d^{3k+2}L(d,k)}{\epsilon}
\]

\noindent satisfies the required condition if the girth of $G$ is
big enough, and all $k$-discs are trees.

\noindent \zlabel{Theorem}{FPS}
\begin{thm}
\noindent (Fichtenberger, Peng and Sohler)

\noindent Let $d\ge2,k\ge1,\epsilon\in(0,1)$. Then for any simple
$d$-bounded graph $G$ with $\text{girth}(G)\ge2k+2$ there is an
simple graph $H$ such that
\[
||\text{freq}_{k}(G)-\text{freq}_{k}(H)||_{1}\le\epsilon\quad\wedge\quad|V(H)|\le36\frac{d^{3k+2}L(d,k)}{\epsilon}
\]
\end{thm}

\noindent 
\noindent In this setting, the $k$-discs of all vertices in $G$
are trees, a fact that allows the authors to utilize the ``Rewire
and Split'' graph manipulation technique to prove the theorem. It
was sketched by the same authors that one can also construct the required
graph $H$ if $G$ is planar by using the planar separator theorem.
\\
\\
A different question with some resemblance to Alon's was posed and
answered by Winkler \zBibRef{W}. It also concerns $k$-discs, but
in the setting of directed edge-colored graphs. Formal definition
of $k$-discs for directed edge-colored graphs is given in chapter
$2$.

\noindent \zlabel{Theorem}{Winkler}
\begin{thm}
\noindent (Winkler)

\noindent There is no deterministic algorithm, that, given $d\ge2,k\ge1$
and a set $\varPhi$ of $d$-bounded $k$-discs of directed, edge-colored
graphs , decides whether there is a directed edge-colored graph $G$
with

\noindent 
\[
\{\text{disc}_{k}(v)|v\in V(G)\}=\varPhi
\]

\noindent In other words, the set of $k$-discs of vertices in $G$
is \textbf{exactly} $\varPhi$.\\
\end{thm}

\noindent 
\noindent The proof by Winkler is based on a reduction from PCP (see
\zgoto{PCP}). A Post Correspondence System (PCS) $P$ is used to
construct a set $\Phi$ of directed edge-colored $k$-discs, such
that there is a graph whose set of $k$-discs is exactly $\varPhi$
if and only if $P$ has a solution. The construction utilizes the
edge directness and coloring to represent the ``letters'' and the
``words'' in $P$. A similar reduction was constructed independently
by Bulitko \zBibRef{BU}, from a slightly different variant of PCP.
It was shown by Jacobs \zBibRef{J} that Winkler's problem is still
undecidable even if $G$ is required to be planar and bipartite.\zSection{Our Contribution}

\noindent Our main result deals with the variant of Alon's question
for directed edge-colored graphs.\zlabel{Question}{DiCo}
\begin{question}
\noindent \textup{(Alon - Directed Edge-Colored Variant)}

\noindent Let $C$ be a finite set of colors, and let $d\ge2,k\ge1,\epsilon\in(0,1)$.
Is there a computable function $f_{C}\coloneqq f_{C}(d,k,\epsilon)$
such that for any $d$-bounded directed graph $G$ whose edges attain
colors in $C$ there is a graph $H$ such that
\[
||\text{freq}_{k}(G)-\text{freq}_{k}(H)||_{1}\le\epsilon\quad\wedge\quad|V(H)|\le f_{C}(d,k,\epsilon)
\]
\end{question}

\noindent We show that this question is interreducible with the original
question for simple graphs.\zlabel{Theorem}{DiCoMain}
\begin{thm}
\noindent Answering \textup{\zgoto{DiCo} }is interreducible with
answering \textup{\zgoto{Main}.}
\end{thm}

This theorem will be a direct corollary of a more general statement,
concerning ``natural'' graph models. \zlabel{Theorem}{InterrInformal}
\begin{thm}
\noindent \textup{(Interreducibility Theorem)}

Given any \textbf{natural} graph model $M$, the variant of Alon's
question for $M$ is interreducible with \textup{\zgoto{Main}}.
\end{thm}

\noindent In Chapter $2$, we will formally define what makes $M$
``natural'' and how the variant of Alon's question is defined.\\

\noindent We will use the tools that we will develop to prove the
Interreducibility Theorem to show that Winkler's question in the simple
graph model is also undecidable.\zlabel{Theorem}{SimpleWinkler}
\begin{thm}
\noindent \textup{(Winkler - Simple Variant)}

\noindent There is no deterministic algorithm, that, given $d\ge2,k\ge1$
and a set $\varPhi$ of $d$-bounded $k$-discs of \textbf{simple}
graphs , decides whether there is a simple graph $G$ with

\noindent 
\[
\{\text{disc}_{k}(v)|v\in V(G)\}=\varPhi
\]

\noindent In other words, the set of $k$-discs of vertices in $G$
is \textbf{exactly} $\varPhi$.\\
\end{thm}

\noindent 
\noindent There is significant resemblance between the questions of
Alon and Winkler. Both questions examine the $k$-disc sets of graphs,
and ask if it is possible to find small graphs that satisfy some restriction
on that set. We therefore conjecture that there is a reduction from
PCP (or a variant of it) to the directed edge-color variant of Alon's
question, and in particular Alon's original question is undecidable.\zlabel{Conjecture}{MainConj}
\begin{conjecture}
\noindent There is no computable function $f(d,k,\epsilon)$ that
satisfies the condition of \textup{\zgoto{Main}}.\\
\end{conjecture}

\noindent In the last chapter, we will examine the variant of the
main question where all the graphs are paths. A solution to a single
PCP system (see \zgoto{PCP}) can be thought of as one long string,
and so it is natural to ask whether this string (i.e. directed edge-labeled
path) can be approximated by a fixed size string. We prove that the
question in this case is decidable.

\noindent \zlabel{Theorem}{PathByPathInformal}
\begin{thm}
\noindent \textup{(Alon - Directed Edge-Colored Path Variant)}

\noindent Let $k\ge1,\epsilon\in(0,1)$ and let $C$ be a finite set
of colors/labels . Let $P$ be a directed path with edge colors in
$C$, then there is a directed edge-colored path $Q$ such that

\noindent 
\[
||\text{freq}_{k}(P)-\text{freq}_{k}(Q)||_{1}\le\epsilon\ \wedge\ |Q|\le24960\frac{8^{k}|S|^{2}(2k)^{6|S|}}{\epsilon^{2}}
\]
\end{thm}

\noindent 
\noindent Moreover, we will show that the problem is still decidable
when alternative definitions for local structure of paths are considered.
For example, the local structure of a vertex in a directed labeled
path can be seen as a single ``string''. In this case the frequency
vector represents the frequency of different ``words'' in the path.
We use \zgoto{PathByPathInformal} to show that the question in this
case is still decidable.\\

\noindent Throughout the thesis we will often use bounds/constants
which are not tight, to improve readability. In general, any value
which only depends on $d,k,\epsilon$ and $|C|$ is considered fixed,
small and negligible in comparison to the graph size $|V|$.

\noindent \zEndChapter

\noindent \noindent \zChapter{S-Graphs and The Interreducbility Theorem}

\noindent The problem of finding small graphs that preserve the local
structure of arbitrarily large graph is not restricted to simple graphs.
It is possible to ask the same question for graphs with one or more
additional properties (directed edges, edge/vertex coloring, multi-edges,
loops and more).

\noindent The motivation behind asking a variant of the question for
other graph types is that adding more ``information'' to the graph
might make it easier to compute the corresponding function $f$ or
prove that it is uncomputable.

\noindent For example, in the scenario of directed graphs, the ``natural''
way to define $k$-discs would be the same as in the simple case,
with the additional requirement that $k$-disc isomorphisms will also
preserve edge direction.\\

\noindent In this chapter, we will formally define what makes a graph
model ``natural''. We will use that definition to formally state
the Interreducibility Theorem (\zgoto{InterrInformal}). Finally,
we will derive \zgoto{DiCoMain}, essentially proving that to show
that Alon's question is undecidable, it is enough to show that the
directed edge-colored version is undecidable.

\zSection{S-Graphs}

When working with non simple graphs (i.e. graphs with some property
like edge coloring), the ``natural'' way to define isomorphism between
two $k$-discs would be by a root preserving isomorphism which also
preserves the property. An important observation here is that there
is nothing special about properties like coloring, directness or multi-edges.
Each such property will only affect the amount of possible $k$-discs,
and not the logic that is used behind their definition. As all properties
will have essentially the same version of the problem, it would be
easier to work with a more general definition and then specify how
each specific property is realized by this definition. To this end,
we introduce the notion of $S$-graphs, as a generalization for graphs
where two $k$-discs are said to be isomorphic if the graph isomorphism
also preserves the additional properties of the model.
\begin{defn}
Let $S$ be a finite non empty set, which we will call the \textbf{information
set}.

Let $V_{S}$ be a finite set of vertices, and let $I$ be a function
\[
I:V_{S}\times V_{S}\to\{0\}\cup\left(\{1\}\times S\right)
\]

We say that the tuple $(V_{S},I)$ is an $S$\textbf{-graph}, and
denote the set of all such tuples by $\Omega(S)$.

For an $S$-graph $G_{S}=(V_{S},I)$, we say that $V_{S}=V_{S}(G_{S})$
is the vertex set of $G_{S}$, and that $I=I(G_{S})$ is the \textbf{information
function }of $G_{S}$.
\end{defn}

The idea behind this definition is that many different graph types/properties
can be defined by choosing the correct information set $S$ and then
defining constrains on the function $I$.

In a sense, the $\{0,1\}$ part of the image of $I$ stands for whether
there is a directed edge from one vertex to another, and the set $S$
contains all the additional information (like edge-coloring, for example).
\begin{example}
(Examples of $S$-graph models)
\end{example}

\begin{itemize}
\item If $S=\{0\}$ then $G_{S}=(V_{S},I)$ can be seen as a directed graph
(where loops are allowed). If we also define that $\forall v\in V\ I(v,v)=0$
then loops are not allowed.
\item If we also require that $\forall v_{1},v_{2}\in V_{S}\ I(v_{1},v_{2})=I(v_{2},v_{1})$,
then every ``edge'' appears in the graph if and only if the reverse
edge appears. In this case the model represents undirected graphs.
\item Edge coloring can be defined by setting $S=\{c_{1},...,c_{m}\}$,
where each element represents a color. In this case, every edge in
the graph will have a single unique color given to it.
\item Edge multiplicity can be defined by taking $S=[t]$ where $t$ is
the maximal edge multiplicity in the graph.
\end{itemize}
In general, any combination of properties can also be represented
by taking suitable $S,I$.\\

Before we can state the variant of Alon's question for $S$-graphs,
we need to go over the basic graph notation. It is important to notice
that most definitions do not depend on $S$, which by itself hints
that the difficulty of the approximation question will not be hindered.
\begin{defn}
Let $S$ be an information set and let $G_{S}=(V_{S},I)$ be an $S$-graph.
\end{defn}

\begin{itemize}
\item Given two distinct vertices $v_{1},v_{2}\in V$, we say that there
is an edge between them if $I(v_{1},v_{2})\ne0$ or $I(v_{2},v_{1})\ne0$.
In this case we say that $v_{1},v_{2}$ are adjacent.
\item The underlying simple graph of $G_{S}$, denoted by $U(G_{S})$ is
defined as the simple graph $G=(V,E)$ that is created by taking $V=V_{S}$
and $E=\{(u,v)|I(u,v)\ne0\}$.
\item The distance between $v_{1},v_{2}$ is the length of the shortest
sequence of edges from $v_{1}$ to $v_{2}$.
\item For an integer $d\in\mathbb{N}$, we say that $G$ has maximal degree
at most $d$ if each vertex in $V$ is part of at most $d$ edges
(a loop edge counts as $2$ edges). In particular, each vertex can
have at most $d$ neighbors.
\item For any $v\in V_{S}$, the $k$-disc of $v$, denoted by $\text{disc}_{k}(v)$
or $\text{disc}_{k}(G_{S},v)$ is defined as the subgraph that is
induced by the vertices that are at distance at most $k$ to $v$
in $G_{S}.$
\item We say that two $k$-discs are\textbf{ }isomorphic if and only if
there is a root-preserving graph isomorphism which also preserves
the function $I$. Namely, two $k$-discs $\Gamma_{1}=(V_{1},I_{1})$
and $\Gamma_{2}=(V_{2},I_{2})$ are isomorphic if and only if there
is a graph isomorphism $f:V_{1}\to V_{2}$ such that 
\[
\forall v_{1},v_{2}\in V_{1}\quad I_{1}(v_{1},v_{2})=I_{2}(f(v_{1}),f(v_{2}))
\]
\end{itemize}
It is important to note that only the definition of $k$-disc isomorphism
depends on $S$. Everything else is exactly the same as in the simple
graph model. It is possible to define $k$-discs differently for some
$S$-graph models; we will examine some alternative definitions in
the last chapter of the thesis.\\
\\
We denote the set of all non isomorphic $d$-bounded degree rooted
$S$-graphs with radius at most $k$ by $\mathcal{L}_{S}(d,k)$. Just
like in the simple case, we have the following fact.

\zlabel{Fact}{SDiscSize}
\begin{fact}
Let $v\in V(G_{S})$ be a vertex of an $S$-graph, then $|\text{disc}_{k}(v)|\le2d^{k}$.
In particular $\mathcal{L_{S}}(d,k)$ is finite.
\end{fact}

The same reasoning as in \zgoto{DiscSize}, together with $S$ being
finite, can be used to prove this fact.\\

We denote the size of $\mathcal{L}_{S}(d,k)$ by $L_{S}\coloneqq L_{S}(d,k)$.
The $k$-disc count and frequency distribution vectors - $\text{cnt}_{k}(G_{S})$
and $\text{freq}_{k}(G_{S})$, are defined in the exact same way as
in the simple case, with the only difference being in the amount of
entries in the vectors ($L_{S}(d,k)$ instead of $L(d,k)$).

\zSection{Properties of S-Graphs}

In this section we will state and prove some very useful lemmas which
will be used as auxiliary properties of $S$-graphs throughout the
rest of the thesis.

We start with a lemma that gives an estimation of the difference between
the frequency distribution of an $S$-graph and one of its subgraphs.

\zlabel{Lemma}{FreqDiff}
\begin{lem}
Let $d\ge2,k\ge1$ and let $S$ be an information set. 

Suppose $G$ is an $S$-graph with maximum degree $d$ and $H$ is
an induced subgraph of $G$. Then
\[
||\text{freq}_{k}(G)-\text{freq}_{k}(H)||_{1}\le\frac{(1+2d^{k})\left(|G|-|H|\right)}{|H|}
\]
\end{lem}

\textbf{Proof }By \zgoto{SDiscSize} we know that removing a single
vertex from $G$ affects the $k$-disc of at most $2d^{k}$ vertices.
In general, removing $x$ vertices from $G$ will affect the $k$-discs
of at most $2d^{k}x$ vertices. In our case, $x=|G|-|H|$ is the amount
of vertices that were removed from $G$, and so for every $k$-disc
$\Gamma\in\mathcal{L}_{S}(d,k)$ it holds that
\[
|\text{cnt}_{k}(G,\Gamma)-\text{cnt}_{k}(H,\Gamma)|\le2d^{k}\left(|G|-|H|\right)
\]

This bound can be generalized to a bound on the frequency distribution
difference

\begin{align*}
 & |G|\cdot|H|\cdot||\text{freq}_{k}(G)-\text{freq}_{k}(H)||_{1}=\\
 & =|G|\cdot|H|\cdot\sum_{\Gamma\in\mathcal{L_{S}}(d,k)}|\text{freq}_{k}(G,\Gamma)-\text{freq}_{k}(H,\Gamma)|=|G|\cdot|H|\cdot\sum_{\Gamma\in\mathcal{L_{S}}(d,k)}|\frac{\text{cnt}_{k}(G,\Gamma)}{|G|}-\frac{\text{cnt}_{k}(H,\Gamma)}{|H|}|=\\
 & =\sum_{\Gamma\in\mathcal{L_{S}}(d,k)}||H|\text{cnt}_{k}(G,\Gamma)-|G|\text{cnt}_{k}(H,\Gamma)|=\sum_{\Gamma\in\mathcal{L_{S}}(d,k)}|\left(|H|-|G|+|G|\right)\text{cnt}_{k}(G,\Gamma)-|G|\text{cnt}_{k}(H,\Gamma)|=\\
 & =\sum_{\Gamma\in L_{S}(d,k)}|\left(|H|-|G|\right)\text{cnt}_{k}(G,\Gamma)+|G|\left(\text{cnt}_{k}(G,\Gamma)-\text{cnt}_{k}(H,\Gamma)\right)|\le\\
 & \le\left(|G|-|H|\right)\cdot\sum_{\Gamma\in\mathcal{L_{S}}(d,k)}\text{cnt}_{k}(G,\Gamma)+|G|\cdot\sum_{\Gamma\in\mathcal{L_{S}}(d,k)}|\text{cnt}_{k}(G,\Gamma)-\text{cnt}_{k}(H,\Gamma)|=\\
 & =\left(|G|-|H|\right)\cdot|G|+|G|\cdot\sum_{\Gamma\in\mathcal{L_{S}}(d,k)}|\text{cnt}_{k}(G,\Gamma)-\text{cnt}_{k}(H,\Gamma)|\le\\
 & \le\left(|G|-|H|\right)\cdot|G|+|G|\cdot2d^{k}\left(|G|-|H|\right)=(1+2d^{k})|G|\left(|G|-|H|\right)
\end{align*}

By isolating the frequency difference we conclude that
\[
||\text{freq}_{k}(G)-\text{freq}_{k}(H)||_{1}\le\frac{(1+2d^{k})|G|\left(|G|-|H|\right)}{|G|\cdot|H|}=\frac{(1+2d^{k})\left(|G|-|H|\right)}{|H|}
\]

This completes the proof of \zgoto{FreqDiff}.\qed\\

In the next lemma, we show that if two $S$-graphs on the same vertex
set are close to each other (i.e. one can be formed by adding/removing/changing
a small amount of edges in the other) then the difference between
their FDVs is small.

\zlabel{Lemma}{EdgeChange}
\begin{lem}
Let $d\ge2,k\ge1$ and let $S$ be an information set.

Suppose $G=(V,I_{G})$ is an $S$-graph with maximum degree $d$ and
$H=(V,I_{H})$ is an $S$-graph formed by adding/removing/changing
$m\ge1$ \textbf{edges} in $G$ (values of the function $I_{G}$).
Then
\[
||\text{freq}_{k}(G)-\text{freq}_{k}(H)||_{1}\le\frac{4d^{k}mL_{S}(d,k)}{|G|}
\]
\end{lem}

\textbf{Proof }Suppose that $H$ was formed by adding/removing/changing
$m$ edges in $G$.

Each affected edge has exactly two end points, and each such end point,
by \zgoto{SDiscSize} , belongs to the $k$-disc of at most $2d^{k}$
vertices. In total the amount of vertices whose $k$-discs have changed
as a result of the single edge change is at most $2\cdot2d^{k}=4d^{k}$.
Therefore, the total amount of affected $k$-discs is at most $4d^{k}m$.
Using the fact that $|G|=|H|$ we have
\begin{align*}
 & ||\text{freq}_{k}(G)-\text{freq}_{k}(H)||_{1}\\
 & =\sum_{\Gamma\in\mathcal{L}_{S}(d,k)}|\text{freq}_{k}(G,\Gamma)-\text{freq}_{k}(H,\Gamma)|=\frac{1}{|G|}\sum_{\Gamma\in\mathcal{L}_{S}(d,k)}|\text{cnt}_{k}(G,\Gamma)-\text{cnt}_{k}(H,\Gamma)|\le\\
 & \le\frac{1}{|G|}\sum_{\Gamma\in\mathcal{L}_{S}(d,k)}4d^{k}m=\frac{4d^{k}mL_{S}(d,k)}{|G|}
\end{align*}

This completes the proof of \zgoto{EdgeChange}.\qed\\

The next lemma is a very powerful tool that will be used throughout
the thesis.

\zlabel{Lemma}{WeightShifting}
\begin{lem}
(Weight Shifting Lemma)

Let $d\ge2,k\ge1$ and let $S$ be an information set. Suppose $G,H_{1},H_{2}$
are $S$-graphs such that for every $k$-disc $\Gamma\in\mathcal{L}_{S}(d,k)$
the following holds
\[
\text{freq}_{k}(H_{2},\Gamma)<\text{freq}_{k}(H_{1},\Gamma)\rightarrow\text{freq}_{k}(G,\Gamma)=0
\]

Then $||\text{freq}_{k}(G)-\text{freq}_{k}(H_{2})||_{1}\le||\text{freq}_{k}(G)-\text{freq}_{k}(H_{1})||_{1}$.
\end{lem}

In other words, if we can create the vector $\text{freq}_{k}(H_{2},\Gamma)$
from $\text{freq}_{k}(H_{1},\Gamma)$ by ``shifting weight'' away
from ``bad'' entries (where $\text{freq}_{k}(G,\Gamma)=0$), then
$H_{2}$ gives a better approximation than $H_{1}$ of the local structure
of $G$.\\
\\
\textbf{Proof }By the definition of the frequency distribution vector
of $G,H_{1},H_{2}$, we have
\begin{equation}
\sum_{\Gamma\in\mathcal{L}_{S}(d,k)}\text{freq}_{k}(G,\Gamma)=\sum_{\Gamma\in\mathcal{L}_{S}(d,k)}\text{freq}_{k}(H_{1},\Gamma)=\sum_{\Gamma\in\mathcal{L}_{S}(d,k)}\text{freq}_{k}(H_{2},\Gamma)=1\label{eq:1}
\end{equation}

We define a partition of $\mathcal{L}_{S}(d,k)$ into two sets
\[
F_{1}=\left\{ \Gamma\in\mathcal{L}_{S}(d,k)|\text{freq}_{k}(H_{2},\Gamma)<\text{freq}_{k}(H_{1},\Gamma)\right\} \quad F_{2}=\mathcal{L}_{S}(d,k)\backslash F_{1}
\]

By the assumption of the lemma and \ref{eq:1}, we have

\begin{align*}
 & ||\text{freq}_{k}(G)-\text{freq}_{k}(H_{2})||_{1}=\\
 & =\sum_{\Gamma\in\mathcal{L}_{S}(d,k)}|\text{freq}_{k}(G,\Gamma)-\text{freq}_{k}(H_{2},\Gamma)|=\\
 & =\sum_{\Gamma\in F_{1}}|\text{freq}_{k}(G,\Gamma)-\text{freq}_{k}(H_{2},\Gamma)|+\sum_{\Gamma\in F_{2}}|\text{freq}_{k}(G,\Gamma)-\text{freq}_{k}(H_{2},\Gamma)|=\\
 & =\sum_{\Gamma\in F_{1}}|0-\text{freq}_{k}(H_{2},\Gamma)|+\sum_{\Gamma\in F_{2}}|\text{freq}_{k}(G,\Gamma)-\text{freq}_{k}(H_{1},\Gamma)+\text{freq}_{k}(H_{1},\Gamma)-\text{freq}_{k}(H_{2},\Gamma)|\le\\
 & \le\sum_{\Gamma\in F_{1}}\text{freq}_{k}(H_{2},\Gamma)+\sum_{\Gamma\in F_{2}}|\text{freq}_{k}(G,\Gamma)-\text{freq}_{k}(H_{1},\Gamma)|+\sum_{\Gamma\in F_{2}}|\text{freq}_{k}(H_{1},\Gamma)-\text{freq}_{k}(H_{2},\Gamma)|=\\
 & =\sum_{\Gamma\in F_{2}}|\text{freq}_{k}(G,\Gamma)-\text{freq}_{k}(H_{1},\Gamma)|+\sum_{\Gamma\in F_{1}}\text{freq}_{k}(H_{2},\Gamma)+\sum_{\Gamma\in F_{2}}\left(\text{freq}_{k}(H_{2},\Gamma)-\text{freq}_{k}(H_{1},\Gamma)\right)=\\
 & =\sum_{\Gamma\in F_{2}}|\text{freq}_{k}(G,\Gamma)-\text{freq}_{k}(H_{1},\Gamma)|+\sum_{\Gamma\in\mathcal{L}_{S}(d,k)}\text{freq}_{k}(H_{2},\Gamma)-\sum_{\Gamma\in F_{2}}\text{freq}_{k}(H_{1},\Gamma)=\\
 & =\sum_{\Gamma\in F_{2}}|\text{freq}_{k}(G,\Gamma)-\text{freq}_{k}(H_{1},\Gamma)|+\sum_{\Gamma\in F_{1}}\text{freq}_{k}(H_{1},\Gamma)=\\
 & =\sum_{\Gamma\in F_{2}}|\text{freq}_{k}(G,\Gamma)-\text{freq}_{k}(H_{1},\Gamma)|+\sum_{\Gamma\in F_{1}}|\text{freq}_{k}(G,\Gamma)-\text{freq}_{k}(H_{1},\Gamma)|=\\
 & =\sum_{\Gamma\in\mathcal{L}_{S}(d,k)}|\text{freq}_{k}(G,\Gamma)-\text{freq}_{k}(H_{1},\Gamma)|\\
 & =||\text{freq}_{k}(G)-\text{freq}_{k}(H_{1})||_{1}
\end{align*}

This completes the proof of \zgoto{WeightShifting}.\qed\\

The last lemma of this section concerns alternative definitions of
local structure of graphs. 

Suppose $M:\mathcal{L}_{S}(d,k)\to X$ is a function that maps the
set of $k$-discs into some finite set $X=\{x_{1},...,x_{|X|}\}$.
Given an $S$-graph $G$, the frequency distribution vector $\text{freq}_{M}(G)$
is the vector where the $i$-th entry counts the fraction of vertices
in $G$ whose $k$-disc attains a value $x_{i}$ by $M$.

\zlabel{Lemma}{FreqDiffModulo}
\begin{lem}
Let $d\ge2,k\ge1$ and let $S$ be an information set. 

Suppose $G,H$ are $S$-graphs, and ${M:\mathcal{L}_{S}(d,k)\to X}$
is a function that maps $k$-discs into some finite set $X$. Then
\[
||\text{freq}_{M}(G)-\text{freq}_{M}(H)||_{1}\le||\text{freq}_{k}(G)-\text{freq}_{k}(H)||_{1}
\]
\end{lem}

In other words, by mapping $\mathcal{L}_{S}(d,k)$ to $X$, we can
only ``lose'' information about the local structure.\\

\textbf{Proof} Let $1\le i\le|X|$. By the definition of $M$, we
have
\[
\text{freq}_{M}(G)_{i}=\sum_{\Gamma\in M^{-1}(x_{i})}\text{freq}_{k}(G,\Gamma)\quad\text{freq}_{M}(H)_{i}=\sum_{\Gamma\in M^{-1}(x_{i})}\text{freq}_{k}(H,\Gamma)
\]

And therefore

\begin{align*}
||\text{freq}_{M}(G)-\text{freq}_{M}(H)||_{1} & =\sum_{i=1}^{|X|}|\text{freq}_{M}(G)_{i}-\text{freq}_{M}(H)_{i}|=\\
 & =\sum_{i=1}^{|X|}|\sum_{\Gamma\in M^{-1}(x_{i})}\text{freq}_{k}(G,\Gamma)-\sum_{\Gamma\in M^{-1}(x_{i})}\text{freq}_{k}(H,\Gamma)|\le\\
 & \le\sum_{i=1}^{|X|}\sum_{\Gamma\in M^{-1}(x_{i})}|\text{freq}_{k}(G,\Gamma)-\text{freq}_{k}(H,\Gamma)|=\\
 & =\sum_{\Gamma\in\mathcal{L}_{S}(d,k)}|\text{freq}_{k}(G,\Gamma)-\text{freq}_{k}(H,\Gamma)|=\\
 & =||\text{freq}_{k}(G)-\text{freq}_{k}(H)||_{1}
\end{align*}

This completes the proof of \zgoto{FreqDiffModulo}. \qed

\newpage
\zSection{The Interreducibility Theorem}

In this section we formally state and prove the Interreducibility
Theorem (\zgoto{InterrInformal}). 

We start by defining the variant of Alon's question for $S$-graphs.

\zlabel{Question}{Ext}
\begin{question}
(Alon - S-Graph Variant)

Let $d\ge2,k\ge1,\epsilon\in(0,1)$ and let $S$ be an information
set. Let $A\subseteq\Omega(S)$ be a set of $d$-bounded $S$-graphs.
Is there a computable function $f_{S,A}\coloneqq f_{S,A}(d,k,\epsilon)$
such that for any $S$-graph $G\in A$ there is an $S$-graph $H\in A$
such that
\[
||\text{freq}_{k}(G)-\text{freq}_{k}(H)||_{1}\le\epsilon\quad\wedge\quad|V(H)|\le f_{S,A}(d,k,\epsilon)
\]
\end{question}

Just like in the simple graph case, it can be shown that there is
a function $f_{S,A}(d,k,\epsilon)<\infty$ which satisfies the condition
of the question (see \zgoto{FiniteCover}), but the proof does not
give any information regarding the size of the approximating graph
$H$. Clearly, for some choices of $S,A$, this function is trivially
computable. For example, if $A\subseteq\Omega(S)$ is finite, then
taking $f_{S,A}(d,k,\epsilon)=\max_{G\in A}|V(G)|$ is sufficient.
We wish to show that for any ``natural'' choice of $S$ and $A$,
the question is interreducible with \zgoto{Main}. We proceed by defining
what make $A$ a natural set.
\begin{defn}
Let $d\ge2,k\ge1$ and let $S$ be an information set. Let $A\subseteq\Omega(S)$
be a set of $d$-bounded $S$-graphs.
\end{defn}

\begin{enumerate}
\item We say that $A$ is \textbf{\uline{natural}}\textbf{ }if for every
$G_{S}\in A$ and $H_{S}\in\Omega(S)$ there exists an $H_{S}^{1}\in A$
with 
\[
||\text{freq}_{k}(G_{S})-\text{freq}_{k}(H_{S}^{1})||_{1}\le||\text{freq}_{k}(G_{S})-\text{freq}_{k}(H_{S})||_{1}\quad\wedge\quad|V(H_{S}^{1})|\le|V(H_{S})|
\]
\item We say that $A$ is a \textbf{\uline{natural extension}} if $A$
is natural and for every $d$-bounded simple graph $G$ there is an
$S$-graph $G_{S}\in A$ with $U(G_{S})=G$ (same underlying simple
graph).\\
\end{enumerate}
We can think of the naturality property as the ``crucial'' property
of any interesting set $A$. If a graph $G_{S}\in A$ is approximated
by some graph $H_{S}\in\Omega(S)$ then we would expect that $H_{S}$
itself is a member of $A$ or very close to being one. For example,
if $A$ is the set of directed graphs without loops, and $H_{S}$
contains loops, then clearly we can remove all those loops and get
an even better approximation of the same size. The concept of a natural
extension ``requires'' a natural set $A$ to contain some variation
of every possible simple graph, essentially making $A$ an infinite
set which is ``at least'' the set of simple graphs. These two properties
contain the ``critical'' difference between simple graphs and other
graph models. If a graph model satisfies these two properties then
we would expect the corresponding approximation problem to be interreducible
with the original one, no matter what the model represents. We can
now formally state the main result of the thesis.

\zlabel{Theorem}{InterrFormal}
\begin{thm}
\textup{(Interreducibility Theorem)}

Let $S$ be an information set, and let $A\subseteq\Omega(S)$ be
a \textbf{natural extension}.

Then, \emph{\zgoto{Main} }and\emph{ \zgoto{Ext} }(for this choice
of $S,A$) are interreducible.
\end{thm}

In other words, if we restate the problem for any naturally defined
graph property or properties, then the difficulty of the question
is not altered. In particular, all natural extensions are interreducible
between each other (via the simple variant). This means that answering
the question for a single natural extension pair $(S,A)$ will answer
the question for any other pair and also the simple case. The proof
of the theorem is a direct corollary of the following two lemmas,
which will be proven in Chapters $3$ and $4$, respectively.

\zlabel{Lemma}{EasyOne}
\begin{lem}
\textup{(Reduction from simple graphs to $S$-graphs)}

Let $d\ge2,k\ge1,\epsilon\in(0,1)$ and let $S$ be an information
set. Let $A\subseteq\Omega(S)$ be a \textbf{natural extension }set
of $S$-graphs. Suppose there is a function $f_{S,A}$ that satisfies
the condition of \textup{\zgoto{Ext},} then $f(d,k,\epsilon)\coloneqq f_{S,A}(d,k,\epsilon)$
satisfies the condition of \emph{\zgoto{Main}.}
\end{lem}

\zlabel{Lemma}{HardOne}
\begin{lem}
\textup{(Reduction from $S$-graphs to simple graphs)}

Let $d\ge2,k\ge1,\epsilon\in(0,1)$ and let $S$ be an information
set. Let $A\subseteq\Omega(S)$ be a \textbf{natural }set of $S$-graphs.
Suppose there is a function $f$ that satisfies the condition of \textup{\zgoto{Main}},
then $f_{S,A}(d,k,\epsilon)\coloneqq f(d_{1},k_{1},\epsilon_{1})$
where
\[
t=\max\{d+4,|S|\}\quad d_{1}=2t+1\quad k_{1}=3k\quad\epsilon_{1}=\frac{\epsilon}{4(2t+2)^{2}(1+2(2t+1)^{q})}
\]

satisfies the condition of\emph{ \zgoto{Ext}.}
\end{lem}

\textbf{Proof of \zgoto{InterrInformal} / \zgoto{InterrFormal}}

Follows immediately from \zgoto{EasyOne} and \zgoto{HardOne}. \qed

\zSection{Examples of Natural Extensions}

The simplest example of a natural extension over an information set
$S$ is the set $A=\Omega(S)$.

\zlabel{Proposition}{WholeSet}
\begin{prop}
Let $S$ be an information set. Then $A=\Omega(S)$ is a natural extension.
\end{prop}

\textbf{Proof }We prove that the two required conditions hold for
\textbf{$A$.}
\begin{enumerate}
\item Let $G_{S}\in A$ and $H_{S}\in\Omega(S)$. Taking $H_{S}^{1}=H_{S}\in\Omega(S)=A$
is sufficient.
\item Let $G$ be a $d$-bounded graph. We need to find an $S$-graph $G_{S}$
with $U(G_{S})=G$. Let $s\in S$ be a value, we then define the graph
$G_{S}=(V(G),I)$ where
\[
I(v_{1},v_{2})=I(v_{2},v_{1})=\begin{cases}
0 & (v_{1},v_{2})\notin E(G)\\
\{1,s\} & (v_{1},v_{2})\in E(G)
\end{cases}
\]
Clearly $U(G_{S})=G$, as required.\\
\end{enumerate}
We observe that $\Omega(S)$ is in fact the set of all directed graphs
(where loops and bidirectional edges are allowed), whose edges are
colored in $|S|$ colors. We now use the proposition to prove \zgoto{DiCoMain}.\\

\textbf{Proof of \zgoto{DiCoMain} }Take the information set $S=C$
and set $A=\Omega(S)$ to be the set of all $S$-graphs. In this setting,
$S$-graphs in $A$ are exactly directed graphs whose edges are colored
by colors in $C$. By \zgoto{WholeSet}, $A$ is a natural extension,
and by the Interreducibility Theorem, the problem for $S,A$ is interreducible
with the simple problem. In the special case $|C|=1$, the problem
is equivalent to asking the question for directed graphs (without
edge colors).\qed\\

Another important family of natural models are graph models where
the definition of the set $A$ only relies on ``local'' restrictions,
like coloring, multiplicity or direction. We formalize this in the
following Lemma.

\zlabel{Lemma}{EdgeNaturals}
\begin{lem}
Let $d\ge2$ and let $S$ be an information set. Let $P\subseteq S\times S$,
$Q\subseteq S$ be some sets such that $(0,0)\in P$ and $0\in Q$.
Suppose that $A\subseteq\Omega(S)$ is the set of all $S$-graphs
$G_{S}=(V_{S},I)$ such that
\[
\bigcup_{u\ne v\in V_{S}}\left(I\left(u,v\right),I\left(v,u\right)\right)\subseteq P\quad\wedge\quad\bigcup_{v\in V_{S}}I\left(v,v\right)\subseteq Q
\]

Then $A$ is \textbf{natural}.
\end{lem}

In other words, if $A$ is the set of \textbf{all} $S$-graphs that
comply to some local restriction on edges between vertices, then $A$
is natural.\\

\textbf{Proof} Let $k\ge1$, let $G_{S}\in A$ be an $S$-graph, and
let $H_{S}\in\Omega(S)$ be an $S$-graph. We need to show that there
is an $S$-graph $H_{S}^{1}\in A$ such that
\[
||\text{freq}_{k}(G_{S})-\text{freq}_{k}(H_{S}^{1})||_{1}\le||\text{freq}_{k}(G_{S})-\text{freq}_{k}(H_{S})||_{1}\quad\wedge\quad|V(H_{S}^{1})|\le|V(H_{S})|
\]

We construct $H_{S}^{1}$ along the execution of the following algorithm

\textbf{\uline{Algorithm}}
\begin{enumerate}
\item \textbf{function} ConstructSGraph($H_{S}=\left(V(H_{S}),I\right)$)
\item \qquad$V_{new}\leftarrow V(H)$
\item \qquad$I_{new}\leftarrow I$
\item \qquad \textbf{for} $(u,v)\in V(H_{S})\times V(H_{S})$ \textbf{do}
\item \qquad\qquad \textbf{if} $u\ne v$ \textbf{then}
\item \qquad\qquad\qquad \textbf{if} $\left(I\left(u,v\right),I\left(v,u\right)\right)\notin P$
\textbf{then}
\item \qquad\qquad\qquad\qquad $I_{new}(u,v)=I_{new}(v,u)=0$
\item \qquad\qquad\qquad \textbf{end if}
\item \qquad\qquad \textbf{else}
\item \qquad\qquad\qquad \textbf{if} $I(v,v)\notin Q$ \textbf{then}
\item \qquad\qquad\qquad\qquad $I_{new}(v,v)=0$
\item \qquad\qquad\qquad \textbf{end if}
\item \qquad\qquad \textbf{end} \textbf{if}
\item \qquad \textbf{end} \textbf{for}
\item \qquad \textbf{return} $H_{S}^{1}\coloneqq(V_{new},I_{new})$
\item \textbf{end} \textbf{function}
\end{enumerate}
The algorithm starts with the set $H_{S}^{1}\coloneqq H_{S}$ and
then goes over all pairs of vertices in $V(H_{S})$. For each such
pair, if the edges between them (or between a vertex and itself) are
not in $P$ or $Q$ respectively, they are removed.

By the definition of $A$, we know that the resulting graph $H_{S}^{1}$
is in $A$, and we also know that $|V(H_{S}^{1})|=|V(H_{S})|$ as
required. It remains to show that
\[
||\text{freq}_{k}(G_{S})-\text{freq}_{k}(H_{S}^{1})||_{1}\le||\text{freq}_{k}(G_{S})-\text{freq}_{k}(H_{S})||_{1}
\]

We prove this by induction on the number of iterations of the for-loop
in the algorithm.

\textbf{Base}: Before the for-loop section of the algorithm, $H_{S}^{1}$
is precisely $H_{S}$ and so the inequality is true.

\textbf{Step}: For readability purposes, denote the $S$-graph before
and after the $n$-th iteration by $H_{S}^{n}$ and $H_{S}^{n+1}$,
respectively. We assume by induction that
\[
||\text{freq}_{k}(G_{S})-\text{freq}_{k}(H_{S}^{n})||_{1}\le||\text{freq}_{k}(G_{S})-\text{freq}_{k}(H_{S})||_{1}
\]

Now, suppose that the $n$-th iteration of the loop considers the
pair $(u,v)\in V(H_{S})\times V(H_{S})$, and suppose that $u\ne v$
(the case $u=v$ is similar).

If $\left(I\left(u,v\right),I\left(v,u\right)\right)\in P$, then
$H_{S}^{n+1}=H_{S}^{n}$ and the so the inequality continues to hold.

If $\left(I\left(u,v\right),I\left(v,u\right)\right)\notin P$, then
any $k$-disc $\text{disc}_{k}(v)$ which was affected the removal
of edges between $u$ and $v$ does not appear in $G_{S}$ (because
$G_{S}\in A$ and so this edge pair cannot appear in it). In other
words, if a $k$-disc $\Gamma$ appears in $H_{S}^{n}$ more than
in $H_{S}^{n+1}$, then it must contain the edges between $u$ and
$v$ and therefore $\text{freq}_{k}(G,\Gamma)=0$. We can thus write
\[
\text{freq}_{k}(H_{S}^{n+1},\Gamma)<\text{freq}_{k}(H_{S}^{n},\Gamma)\rightarrow\text{freq}_{k}(G,\Gamma)=0
\]

This is exactly the required condition in the weight-shifting lemma
(\zgoto{WeightShifting}), and therefore
\[
||\text{freq}_{k}(G_{S})-\text{freq}_{k}(H_{S}^{n+1})||_{1}\le||\text{freq}_{k}(G_{S})-\text{freq}_{k}(H_{S}^{n})||_{1}
\]

Which is what we had to prove. This completes the induction and the
proof of \zgoto{EdgeNaturals}.\qed\\

We give some concrete examples of sets $P,Q$ that satisfy the conditions
of the lemma.

\zlabel{Example}{NatExtExample}
\begin{example}
(Examples of graph models that can be realized by appropriate $S,P,Q$)
\end{example}

\begin{enumerate}
\item \textbf{\uline{Simple graphs}}\\
Take $S=\{0\}$ and define $P=\{\left(0,0\right),\left(\left(1,0\right),\left(1,0\right)\right)\},Q=\{0\}$.
The set $A$ is then exactly $S$-graphs such that for every vertex
pair, there are either no edges at all or a pair of directional edges
attaining the value $0$. We can create an isomorphism between each
such graph with the underlying simple graph, where each edge pair
is replaced with an undirected edge. The underlying simple graph has
no loops as we have chosen $Q=\{0\}$. In particular, anything we
prove for $S$-graphs is also true for simple graphs.
\item \textbf{\uline{Directed graphs (no loops, no bidirectional edges)
with colored edges}}\\
Given a set $C$ of edge colors, take $S=C$ and $P=\left(0,0\right)\cup\bigcup_{c\in C}\{\left(\left(1,c\right),0\right),\left(0,\left(1,c\right)\right)\},Q=\{0\}$.
Loops cannot occur by the definition of $Q$, and bidirectional edges
cannot occur by the definition of $P$.
\item \textbf{\uline{Directed graphs with multiple directed edges between
vertices}}\\
Suppose we allow at most $k$ directed edges from a vertex to another.
We take $S=[k]$ and define the set $M=\{0\}\cup\bigcup_{1\le i\le k}\{(1,i)\}$.
Finally setting $P=M\times M,Q=M$ realizes the required model.\\
\end{enumerate}
All of these models satisfy the second condition of natural extensions
(any simple graph can be represented by a directed/colored graph).
By \zgoto{EdgeNaturals} the corresponding sets $A$ are natural,
and therefore by the Interreducibility Theorem the corresponding questions
for these models are interreducible with \zgoto{Main}.\\

Many more models can be shown to be natural extensions. For example,
vertex coloring is also a natural extension but it requires a more
complicated version of \zgoto{EdgeNaturals}. In general, any combination
of properties that we have shown to be natural extensions is also
a natural extension.

\noindent \zEndChapter

\noindent \noindent \zChapter{Reduction From Simple Graphs to S-Graphs}

\noindent In this chapter we prove \zgoto{EasyOne}, which is the
easier of the two reductions needed by the Interreducibility Theorem.
The lemma states that a function $f_{S,A}$ which satisfies the condition
of \zgoto{Ext} can be used to construct a function $f$ that satisfies
the condition of \zgoto{Main}.

\noindent The main idea that we will use to prove the lemma is that
any simple graph can be ``embedded'' in an $S$-graph with the same
underlying graph structure. This is essentially the second criteria
of the definition of natural extensions. We will need the following
lemma, in which we prove the frequency difference between two $S$-graphs
is at least as big as the frequency difference between their underlying
simple graphs by a direct application of \zgoto{FreqDiffModulo}.

\noindent \zlabel{Lemma}{FreqDiffEasy}
\begin{lem}
\noindent Let $d\ge2,k\ge1$ and let $S$ be an information set. Suppose
$G_{S},H_{S}$ are $S$-graphs, then
\[
||\text{freq}_{k}(U(G_{S}))-\text{freq}_{k}(U(H_{S}))||_{1}\le||\text{freq}_{k}(G_{S})-\text{freq}_{k}(H_{S})||_{1}
\]
\end{lem}

\noindent \textbf{Proof} We look at the function $f:\mathcal{L}_{S}(d,k)\to\mathcal{L}(d,k)$
where $f(\Gamma)=U(\Gamma)$. By \zgoto{FreqDiffModulo} we have

\noindent 
\begin{equation}
||\text{freq}_{f}(G_{S})-\text{freq}_{f}(H_{S})||_{1}\le||\text{freq}_{k}(G_{S})-\text{freq}_{k}(H_{S})||_{1}\label{eq:2}
\end{equation}

\noindent By \zgoto{DiscSize} we know that $\mathcal{L}(d,k)$ is
finite, so we can write $\mathcal{L}(d,k)=\{\Gamma_{1},...,\Gamma_{L}\}$.
We claim that 
\begin{equation}
\text{freq}_{f}(G_{S})=\text{freq}_{k}(U(G_{S}))\label{eq:3}
\end{equation}
It is enough to prove that the equality holds for each coordinate.
Indeed, by the definition of $U(G_{S})$, we have

\noindent 
\[
\forall i\quad\text{freq}_{k}(U(G_{S}))_{i}=\sum_{\Gamma\in\mathcal{L}_{S}(d,k)\wedge U(\Gamma)=\Gamma_{i}}\text{freq}_{k}(G_{S},\Gamma)=\sum_{\Gamma\in f^{-1}(\Gamma_{i})}\text{freq}_{k}(G_{S},\Gamma)=\text{freq}_{f}(G_{S})_{i}
\]

\noindent The same equality holds for $H_{S}$. By \ref{eq:2} and
\ref{eq:3}, we conclude that
\[
||\text{freq}_{k}(U(G_{S}))-\text{freq}_{k}(U(H_{S}))||_{1}=||\text{freq}_{f}(G_{S})-\text{freq}_{f}(H_{S})||_{1}\le||\text{freq}_{k}(G_{S})-\text{freq}_{k}(H_{S})||_{1}
\]

\noindent Which completes the proof of \zgoto{FreqDiffEasy}.\qed\\

\noindent The proof of \zgoto{EasyOne} is a direct application of
the above lemma.\\

\noindent \textbf{Proof of \zgoto{EasyOne}. }Let $d\ge2,k\ge1,\epsilon\in(0,1)$.
Let $S$ be an information set and let $A\subseteq\Omega(S)$ be a
natural extension set. We know that there is a function $f_{S,A}(d,k,\epsilon)$
which satisfies the required condition of \zgoto{Ext}. We claim that
$f(d,k,\epsilon)\coloneqq f_{S,A}(d,k,\epsilon)$ satisfies the condition
in \zgoto{Main}.

\noindent Let $G=(V,E)$ be a simple graph. Our goal is to prove that
there exists a graph with at most $f(d,k,\epsilon)$ vertices that
preserves the local structure of $G$.

\noindent By the definition of a natural extension, we know that there
is an $S$-graph $G_{S}\in A$ such that 
\[
U(G_{S})=G
\]

\noindent (the underlying simple graph of $G_{S}$ is $G$).

\noindent By the definition of $f_{S,A}(d,k,\epsilon)$, there exists
an $S$-graph $H_{S}\in\Omega(S)$ with
\[
||\text{freq}_{k}(G_{S})-\text{freq}_{k}(H_{S})||_{1}\le\epsilon\quad\wedge\quad|V(H_{S})|\le f_{S,A}(d,k,\epsilon)
\]

\noindent We have found a small graph $H_{S}$ that approximates $G_{S}$,
but we do not necessarily know that $H_{S}\in A$. However, by the
naturality property of $A$, we know that there exists a $S$-graph
$H_{S}^{1}\in A$ with 
\[
||\text{freq}_{k}(G_{S})-\text{freq}_{k}(H_{S}^{1})||_{1}\le||\text{freq}_{k}(G_{S})-\text{freq}_{k}(H_{S})||_{1}\wedge|V(H_{S}^{1})|\le|V(H_{S})|
\]

\noindent Denote the underlying graph $U(H_{S}^{1})$ by $H$. This
is a simple graph on $|V(H_{S}^{1})|$ vertices with maximum degree
$d$. By \zgoto{FreqDiffEasy} we have
\begin{align*}
||\text{freq}_{k}(G)-\text{freq}_{k}(H)||_{1} & =||\text{freq}_{k}(U(G_{S}))-\text{freq}_{k}(U(H_{S}^{1}))||_{1}\le\\
 & \le||\text{freq}_{k}(G_{S})-\text{freq}_{k}(H_{S}^{1})||_{1}\le||\text{freq}_{k}(G_{S})-\text{freq}_{k}(H_{S})||_{1}\le\epsilon
\end{align*}

\noindent We also know that $H$ is small
\[
|V(H)|=|V(H_{S}^{1})|\le|V(H_{S})|\le f_{S,A}(d,k,\epsilon)
\]

\noindent And so for an arbitrary $G$ we have found a small graph
$H$ with at most $f(d,k,\epsilon)=f_{S,A}(d,k,\epsilon)$ vertices
which approximates its local structure, which is what we had to prove.\qed

\noindent \zEndChapter

\noindent \noindent \zChapter{Representing S-Graphs by Simple Graphs}

\noindent The main goal of this chapter is to prove \zgoto{HardOne},
which is the harder of the two reductions used to prove the Interreducibility
Theorem. The lemma states that a function $f$ that satisfies the
condition of \zgoto{Main} can be used to construct a function $f_{S,A}$
that satisfies the condition of \zgoto{Ext}.

\noindent The main idea that will be used to prove the lemma is that
it is possible to represent $S$-graphs by simple graphs in a way
that allows reconstructing the original $S$-graph and also somewhat
preserving its local structure. We will construct a transformation
between $S$-graphs and simple graphs that converts each vertex into
an ``ordered cluster'' with undirected edges that preserves the
information of the $S$-graph.\\

\noindent For the rest of this chapter, we will be working with the
information set $S=\{s_{1},...,s_{|S|}\}$ and some natural set $A\subseteq\Omega(S)$
of $S$-graphs with maximum degree $d$, where we will be examining
$k$-discs. In addition, we introduce the following two parameters
which will be used throughout the chapter:

\noindent 
\[
t\coloneqq t(d,S)\coloneqq\max\{\lceil\frac{d}{2}\rceil+3,|S|+1\}\quad q\coloneqq q(k)\coloneqq3k+1
\]

\noindent \zSection{The Transformation $T_{S}$}
\begin{defn}
For every $S$-graph $G_{S}=(V_{S},I)$, we define the transformation
$G\coloneqq T_{S}(G_{S})$ by
\[
G\coloneqq T_{S}(G_{S})\coloneqq(V,E)
\]

Which is a \textbf{simple} graph constructed in the following way:
\end{defn}

\begin{itemize}
\item For every vertex $v\in V_{S}$, we define a cluster of $2t+2$ unique
vertices, denoted by 
\[
\text{cluster}(v)\coloneqq\{v_{in}^{1},...,v_{in}^{t},v_{out}^{1},...,v_{out}^{t},v_{center},v_{marker}\}\coloneqq\{v^{1},...,v^{2t},v_{c},v_{m}\}
\]
\item The vertex set $V$ of $G$ is then defined as the following disjoint
union:
\[
V\coloneqq\dot{\bigcup_{v\in V_{S}}}\text{cluster}(v)
\]
\item For every $v\in V_{S}$, the following edges are added inside $\text{cluster}(v)$:
\begin{itemize}
\item An edge between $v_{c}$ and $v^{i}$ for all $1\le i\le2t$
\item An edge between $v_{m}$ and $v^{i}$ for all $1\le i\le2t-1$
\item An edge between $v^{i}$ and $v^{i+1}$ for all $1\le i\le2t-1$
\item An edge between $v_{c}$ and $v_{m}$
\end{itemize}
\item Given $v,w\in V_{S}$ (not necessarily $v\ne w$), if $I(v,w)=\{1,s_{i}\}$
then we add the edge $(v_{out}^{i},w_{in}^{i})$.\\
\end{itemize}
\begin{figure}[h]
\includegraphics[scale=0.9]{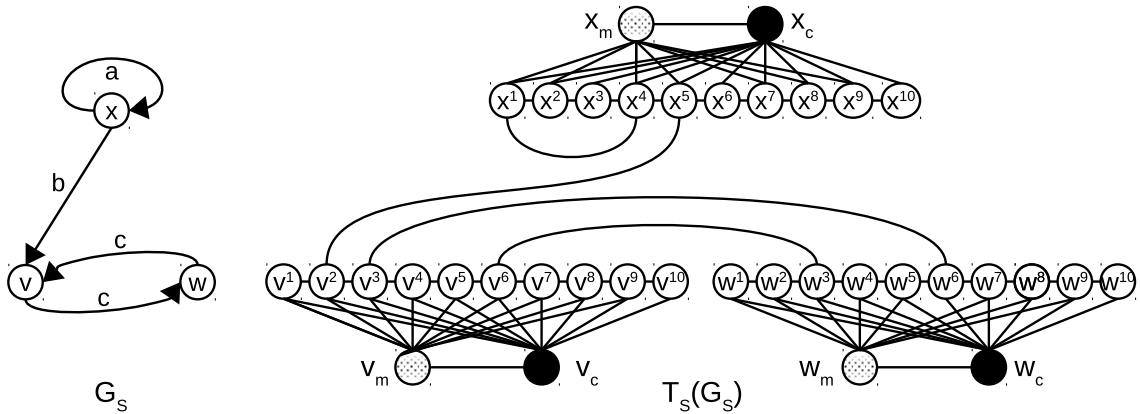}\label{Figure_Ts_Example}

\caption{An example of $G_{S}$ and $T_{S}(G_{S})$ for $d=3,S=\{a,b,c\}$
and $t=5$.}
\end{figure}

The idea behind the construction is that every vertex $v$ in $G_{S}$
is represented by a unique fixed size cluster where each member ``plays''
a very specific role. Each cluster has the exact same structure, containing
four types of vertices: incoming, outgoing, center and marker. 

The incoming/outgoing vertices allow representing edges in an $S$-graph
that attain some value. For example, an edge between $v_{i}^{out}$
and $w_{i}^{in}$ represents an edge in $G_{S}$ between $v$ and
$w$ that attains the value $I(v,w)=s_{i}$. The center vertex is
used to identify which vertices belong to a single cluster, and the
marker vertex is used to distinguish between $v^{1}$ and $v^{2t}$,
which together with the remaining edges instills a unique order in
the cluster. We denote the set of all center-type and marker-type
vertices by $V_{c}$ and $V_{m}$, respectively. The set of all incoming
and outgoing vertices for $s_{i}$ are denoted by $V_{in}^{i}$ and
$V_{out}^{i}$, respectively. We can then write

\[
V=V_{c}\cup V_{m}\cup\bigcup_{1\le i\le t}V_{in}^{i}\cup\bigcup_{1\le i\le t}V_{out}^{i}=V_{c}\cup V_{m}\cup\bigcup_{1\le i\le2t}V^{i}
\]

We claim that $T_{S}$ is well defined, and give some of its basic
properties in the following lemma.

\zlabel{Lemma}{TsProp}
\begin{lem}
Let $G\coloneqq T_{S}(G_{S})\coloneqq(V,E)$ be the transformation
graph of $G_{S}=(V_{S},I)$, then
\end{lem}

\begin{enumerate}
\item $G$ is well defined.
\item $\forall i\quad|V^{i}|=|V_{c}|=|V_{m}|=|V_{S}|=\frac{1}{2t+2}|V|$
\item Let $v\in V_{S}$. Then $\deg_{G}(v_{c})=2t+1$, $\deg_{G}(v_{m})=2t$
and $\forall i\quad\deg_{G}(v^{i})<2t$. In particular the maximal
degree of $G$ is exactly $2t+1$.
\item If $v,w\in V$ have the same $2$-disc, then they are of the same
vertex type ($V^{i}/V_{c}/V_{m}$)
\end{enumerate}
\textbf{Proof}
\begin{enumerate}
\item By definition, the vertex set $V$ is well defined, containing exactly
$|V_{S}|$ disjoint clusters of vertices, each of size $2t+2$. The
only part of the construction which needs careful observation is the
definition of inter-cluster edges. Given an edge $I(v,w)=\{1,s_{i}\}$
in $G_{S}$, we wish to add the edge $(v_{out}^{i},w_{in}^{i})$.
This is only valid if $i\le t$, which is of course true as $i\le|S|<\max\{\lceil\frac{d}{2}\rceil+3,|S|+1\}=t$.
In the special case $v=w$, edges from a vertex to itself correspond
to the edge $(v_{out}^{i},v_{in}^{i})$. We also have the edge $(v_{in}^{t},v_{out}^{1})$,
which is not an inter-cluster edge, but $|S|<t$ and so there can't
be inter-cluster edges that connect to $v_{in}^{t}$.
\item Each vertex in $V_{S}$ is converted to unique $2t+2$ vertices and
so $|V_{S}|=\frac{1}{2t+2}|V|$. Each cluster has exactly $2t+2$
vertex types and therefore $\forall i\quad|V^{i}|=|V_{c}|=|V_{m}|=\frac{1}{2t+2}|V|$.
\item Let $v\in V_{S}$ be a vertex. We examine each case separately:
\begin{itemize}
\item The center and marker vertices $v_{c}$ and $v_{m}$ are part of a
fixed amount of edges inside the cluster, and therefore ${\deg_{G}(v_{c})=2t+1}$,
$\deg_{G}(v_{m})=2t$.
\item If we look at $v_{in}^{i}$ ($1\le i\le t$), then there are at most
$d$ edges between $v_{in}^{i}$ and vertices in other clusters (as
the maximum amount of incoming edges of $v$ in $G_{S}$ is $d$).
Inside the cluster, $v_{in}^{i}$ is connected to at most $4$ vertices:
$v_{c},v_{m},v_{in}^{i-1},v_{in}^{i+1}$. Overall:
\[
\deg(v_{in}^{i})\le d+4=2\left(\frac{d}{2}+2\right)\le2\left(\left(t-3\right)+2\right)=2t-2<2t
\]
\item The same reasoning works for $v_{out}^{i}$ ($1\le i\le t$).
\end{itemize}
\item Suppose $w,x\in V$ are two different vertices with the same $2$-disc.
In particular, $\deg(w)=\deg(x)$. By the previous item, if $\deg(w)=2t+1$
or $\deg(w)=2t$ then both vertices are in $V_{c}$ or $V_{m}$, accordingly.
Otherwise, we know that there are $i,j$ such that $w$ and $x$ are
in $V^{i}$ and in $V^{j}$, respectively. \\
If we look at the $2$-disc of $w$, then it must have exactly one
neighbor with degree $2t+1$, which is the center vertex of the cluster
of $w$ (an inter cluster edge will connect $w$ to a vertex in a
different cluster whose degree is less than $2t$). Moreover, $w$
has at most one neighbor with degree $2t$, which is the marker vertex
of its cluster. The center and marker vertices together define exactly
the order of the other $2t$ vertices in the cluster (the marker distinguishes
the vertex $v^{1}$ from $v^{2t}$ and then each $v^{i}$ defines
$v^{i+1}$). Overall, the entire cluster of $w$ is in it's $2$-disc
and so the path $w\to v_{c}\to v^{1}\to v^{2}\to...\to w$ uniquely
defines the location of $w$ in the cluster. As the $2$-discs of
$w,x$ are the same, this path must be the same, and so $i=j$.\\
\\
\end{enumerate}
Next, we claim that the distance between two center-type vertices
in the transformation graph is at least $3$. In other words, their
$1$-discs do not intersect.

\zlabel{Lemma}{DistThree}
\begin{lem}
Let $G_{S}=(V_{S},I)$ be an $S$-graph, and let $v,w\in V_{S}$ be
two different vertices. Denote the graph $T_{S}(G_{S})$ by $G$.
Then
\begin{enumerate}
\item If $v,w$ are adjacent then $\text{dist}_{G}(v_{c},w_{c})=3$
\item If $v,w$ are not adjacent then $\text{dist}_{G}(v_{c},w_{c})>3$
\end{enumerate}
In particular, $\text{disc}_{1}(v_{c})\cap\text{disc}_{1}(w_{c})=\emptyset$.
\end{lem}

\textbf{Proof }The $1$-discs of $v_{c}$ and $w_{c}$ in $G$ are
both clusters of $2t+2$ vertices, which are disjoint by definition.
The only edges between the clusters are those that connect in-type
and out-type vertices (and not the center-type vertices), and therefore
$\text{dist}_{G}(v_{c},w_{c})\ge3$. Equality is reached only if an
in-type vertex of $v$ and an out-type vertex of $w$ are connected
(or the opposite), which is equivalent to saying that $v,w$ are adjacent.\\

Having established basic properties of the transformation, we want
to connect the local structures of $G_{S}$ and $G$. We know that
each edge $(v,w)$ with value $s$ in $G_{S}$ corresponds to the
path $v_{c}\to v_{out}^{s}\to w_{in}^{s}\to w_{c}$ in $G$, and so
we might hope that the $3k$-discs of center-type vertices in $G$
might be similar to their $k$-disc counterparts in $G_{S}$.

\zSection{The Projection Set $P_q(\Gamma_S)$}

Suppose $v,w\in V_{S}$ are vertices with the same $k$-disc $\Gamma_{S}=\text{disc}_{k}(v)=\text{disc}_{k}(w)$.
It is not necessarily true that $\text{disc}_{3k}(v_{c})=\text{disc}_{3k}(w_{c})$,
as the $3k$-discs in the transformation graph may contain more than
just the center-type vertices that correspond to the $k$-discs in
$G$. However, the $3k$-discs are essentially the same if we think
about the underlying $k$-discs that they represent. To formalize
this idea, we define the projection set of a $k$-disc.
\begin{defn}
Let $\Gamma_{S}\in\mathcal{L}_{S}(d,k)$. The \textbf{$q$}-\textbf{projection
set} $P_{q}(\Gamma_{S})$ of $\Gamma_{S}$ is defined by
\[
P_{q}(\Gamma_{S})\coloneqq\bigcup_{G_{S}=(V_{S},I)\in A}\left\{ \text{disc}_{q}(v_{c})|v\in V_{S},\quad\text{disc}_{k}(v)=\Gamma_{S}\right\} 
\]

We denote by $P_{q}\left(\mathcal{L}_{S}(d,k)\right)$ the set of
all $q$-projections

\[
P_{q}\left(\mathcal{L}_{S}(d,k)\right)\coloneqq\bigcup_{\Gamma_{S}\in\mathcal{L}_{S}(d,k)}P_{q}(\Gamma_{S})
\]
\end{defn}

By \zgoto{TsProp} $(3)$ we know that the maximum degree of transformation
graphs is exactly $2t+1$ and therefore $P_{q}(\Gamma_{S}),P_{q}\left(\mathcal{L}_{S}(d,k)\right)\subseteq\mathcal{L}(2t+1,q)$.
The set $\mathcal{L}(2t+1,q)$ is finite (\zgoto{DiscSize}) and so
$P_{q}(\Gamma_{S}),P_{q}\left(\mathcal{L}_{S}(d,k)\right)$ are finite
as well. We claim that $P_{q}\left(\mathcal{L}_{S}(d,k)\right)$ is
in fact a disjoint union.

\zlabel{Lemma}{ProjDisjointUnion}
\begin{lem}
Let $\Gamma\in P_{q}\left(\mathcal{L}_{S}(d,k)\right)$, then there
is \textbf{exactly one} $\Gamma_{S}$ such that $\Gamma\in P_{q}(\Gamma_{S})$.
\end{lem}

In particular $P_{q}\left(\mathcal{L}_{S}(d,k)\right)$ can be written
as the disjoint union
\[
P_{q}\left(\mathcal{L}_{S}(d,k)\right)=\dot{\bigcup_{\Gamma_{S}\in\mathcal{L}_{S}(d,k)}}P_{q}(\Gamma_{S})
\]

\textbf{Proof} Let $\Gamma\in P_{q}\left(\mathcal{L}_{S}(d,k)\right)$
be a $q$-disc. By definition, there is an $S$-graph $G_{S}=(V_{S},I)$
and a vertex $v\in V_{S}$ such that $\text{disc}_{q}(v_{c})=\Gamma$
in $G\coloneqq T_{S}(G_{S})$. We can ``reconstruct'' $\text{disc}_{k}(v)$
from $\text{disc}_{q}(v_{c})$ with the following deterministic algorithm:

\textbf{\uline{Algorithm}}
\begin{enumerate}
\item \textbf{function} ReconstructKDisc($\Gamma=(V_{\Gamma},E_{\Gamma})$)
\item \qquad$V_{disc}\leftarrow\{v\}$
\item \qquad Init the information function $I_{disc}$ on $V_{disc}$ with
$I_{disc}(v,v)=0$
\item \qquad \textbf{for} $i$ \textbf{in $[k]$ do}
\item \qquad\qquad \textbf{$V_{disc}^{i}=V_{disc}$}
\item \qquad\qquad \textbf{for $w,z$ in $V_{disc}\times V_{\Gamma}$
do}
\item \qquad\qquad\qquad \textbf{if} $\text{dist}(w_{c},z_{c})=3$ \textbf{then}
\item \qquad\qquad\qquad\qquad $V_{disc}^{i}=V_{disc}^{i}\cup\{z\}$
\item \qquad\qquad\qquad \textbf{end if}
\item \qquad\qquad \textbf{end} \textbf{for}
\item \qquad\qquad \textbf{$V_{disc}=V_{disc}^{i}$}
\item \qquad\qquad \textbf{for $w,z,j$ in $V_{disc}\times V_{disc}\times[t]$
do}
\item \qquad\qquad\qquad \textbf{if} $(w_{c},w_{out}^{j}),(w_{out}^{j},z_{in}^{j}),(z_{in}^{j},z_{c})\in E_{\Gamma}$
\textbf{then}
\item \qquad\qquad\qquad\qquad $I_{disc}(w,z)=s_{j}$
\item \qquad\qquad\qquad \textbf{end if}
\item \qquad\qquad \textbf{end} \textbf{for}
\item \qquad \textbf{end} \textbf{for}
\item \qquad \textbf{return} $G_{disc}\coloneqq(V_{disc},I_{disc})$
\item \textbf{end} \textbf{function}
\end{enumerate}
The algorithm starts with only the vertex $v$. In each iteration
of the main for-loop (line $4$), the algorithm performs two steps:
\begin{enumerate}
\item Add any vertex $z\in V_{S}$ such that $\text{dist}(w_{c},z_{c})=3$
for some $w\in V_{S}$. By \zgoto{DistThree}, we know that this condition
is equivalent to saying that $w,z$ are adjacent in $G_{S}$. In other
words, this part of the algorithm adds all vertices in $V_{S}$ that
are adjacent to vertices in $V_{disc}$.
\item Set $I_{disc}(w,z)=s_{j}$ for any pair of vertices $w,z\in V_{disc}$
with $(w_{c},w_{out}^{j}),(w_{out}^{j},z_{in}^{j}),(z_{in}^{j},z_{c})\in E_{\Gamma}$.
In other words, this part of the algorithm adds all the edges between
vertices in $V_{disc}$ from the original graph $G_{S}$.
\end{enumerate}
Overall, after the $i$-th iteration of the algorithm, $(V_{disc},I_{disc})$
is exactly the $i$-disc of $v$. In particular after all $k$ iterations
the algorithm returns the entire $k$-disc.\\
\\
We observe that the algorithm does not depend on the graph $G$ that
we have chosen (we only use $G$ to prove the correctness of the algorithm),
and therefore the $k$-disc of $v$ is uniquely derived from $\Gamma$.
In other words, there is exactly one $\Gamma_{S}$ (the output of
the algorithm) such that $\Gamma\in P_{q}(\Gamma_{S})$, as required.
\qed\\

One can think of the projection set as the set of all $q$-discs in
$\mathcal{L}(2t+1,q)$ that represent the same $k$-disc in $G_{S}$,
enclosed by a single definition. In the next lemma we show how this
definition provides a way to connect the local structures of an $S$-graph
and its transformation.

\zlabel{Lemma}{ProjProp}
\begin{lem}
Let $G_{S}=(V_{S},I)$ be an $S$-graph and let $G\coloneqq T_{S}(G_{S})\coloneqq(V,E)$
be the transformation graph. Then
\end{lem}

\begin{enumerate}
\item Let $w\in V$. Then
\[
\text{disc}_{q}(w)\in P_{q}\left(\mathcal{L}_{S}(d,k)\right)\iff w\in V_{c}
\]
\item The sum of frequencies of $q$-discs in the projection set is given
by
\[
\sum_{\Gamma\in P_{q}\left(\mathcal{L}_{S}(d,k)\right)}\text{freq}_{q}(G,\Gamma)=\frac{1}{2t+2}
\]
\item Given $\Gamma_{S}\in\mathcal{L}_{S}(d,k)$ and a vertex $v\in V_{S}$,
\[
\text{disc}_{k}(v)=\Gamma_{S}\iff\text{disc}_{q}(v_{c})\in P_{q}(\Gamma_{S})
\]
\item Given $\Gamma_{S}\in\mathcal{L}_{S}(d,k)$, the following holds for
the counting vectors
\[
\text{cnt}_{k}(G_{S},\Gamma_{S})=\sum_{\Gamma\in P_{q}(\Gamma_{S})}\text{cnt}_{q}(G,\Gamma)
\]
\item Given $\Gamma_{S}\in\mathcal{L}_{S}(d,k)$, the following holds for
the frequency distribution vectors
\[
\text{freq}_{k}(G_{S},\Gamma_{S})=(2t+2)\cdot\sum_{\Gamma\in P_{q}(\Gamma_{S})}\text{freq}_{q}(G,\Gamma)
\]
\item Let $H_{S}\in\Omega(S)$ be an $S$-graph, and let $H\coloneqq T_{S}(H_{S})$
be its transformation graph, then
\[
||\text{freq}_{k}(H_{S})-\text{freq}_{k}(G_{S})||_{1}\le(2t+2)\cdot||\text{freq}_{q}(H)-\text{freq}_{q}(G)||_{1}
\]
\end{enumerate}
\textbf{Proof}
\begin{enumerate}
\item By definition, $P_{q}\left(\mathcal{L}_{S}(d,k)\right)$ contains
all possible $q$-discs of center-type vertices. If $\text{disc}_{q}(w)\in P_{q}\left(\mathcal{L}_{S}(d,k)\right)$
then it has the same $2$-disc ($2\le q$) as some center-type vertex,
and by \zgoto{TsProp} $(4)$ we know that $w$ is center-type itself,
thus $w\in V_{c}$.
\item By the previous item, we know that any vertex with a $q$-disc in
$P_{q}\left(\mathcal{L}_{S}(d,k)\right)$ is a center-type vertex.
Therefore
\begin{align*}
\sum_{\Gamma\in P_{q}\left(\mathcal{L}_{S}(d,k)\right)}\text{freq}_{q}(G,\Gamma) & =\frac{1}{|G|}\sum_{\Gamma\in P_{q}\left(\mathcal{L}_{S}(d,k)\right)}\text{cnt}_{q}(G,\Gamma)=\frac{1}{|G|}\sum_{\Gamma\in P_{q}\left(\mathcal{L}_{S}(d,k)\right)}|\{v\in V|\text{disc}_{q}(v)=\Gamma\}|\\
 & =\frac{1}{|G|}\cdot|\{v_{c}|v_{c}\in V\}|=\frac{1}{|G|}\cdot|V_{c}|=
\end{align*}
Using \zgoto{TsProp} $(2)$, we get
\[
=\frac{1}{|G|}\cdot\frac{|G|}{2t+2}=\frac{1}{2t+2}
\]
\item The key observation here is that $\text{disc}_{q}(v_{c})\in P(\text{disc}_{k}(v))$
(by the definition of the projection set). Thus:
\begin{itemize}
\item If $\text{disc}_{k}(v)=\Gamma_{S}$, then $\text{disc}_{q}(v_{c})\in P(\text{disc}_{k}(v))=P_{q}(\Gamma_{S})$
\item If $\text{disc}_{q}(v_{c})\in P_{q}(\Gamma_{S})$, then $\text{disc}_{q}(v_{c})\in P(\text{disc}_{k}(v))\cap P_{q}(\Gamma_{S})$.
By \zgoto{ProjDisjointUnion} we get $\text{disc}_{k}(v)=\Gamma_{S}$.
\end{itemize}
\item By the previous item, we have
\begin{align*}
\text{cnt}_{k}(G_{S},\Gamma_{S}) & =|\{v\in V_{S}|\text{disc}_{k}(v)=\Gamma_{S}\}|=|\{v\in V_{S}|\text{disc}_{q}(v_{c})\in P_{q}(\Gamma_{S})\}|=\\
 & =|\{v_{c}\in V_{c}|\text{disc}_{q}(v_{c})\in P_{q}(\Gamma_{S})\}|=
\end{align*}
By the first item in this lemma, we know that only center-type vertices
can have $q$-discs in $P_{q}\left(\mathcal{L}_{S}(d,k)\right)$,
and therefore
\[
=|\{v\in V|\text{disc}_{q}(v_{c})\in P_{q}(\Gamma_{S})\}|=\sum_{\Gamma\in P_{q}(\Gamma_{S})}\text{cnt}_{q}(G,\Gamma)
\]
\item Using the previous item and \zgoto{TsProp} $(2)$, we get
\[
\text{freq}_{k}(G_{S},\Gamma_{S})=\frac{\text{cnt}_{k}(G_{S},\Gamma_{S})}{|V_{S}|}=\frac{\sum_{\Gamma\in P_{q}(\Gamma_{S})}\text{cnt}_{q}(G,\Gamma)}{\frac{|V|}{2t+2}}=(2t+2)\cdot\sum_{\Gamma\in P_{q}(\Gamma_{S})}\text{freq}_{q}(G,\Gamma)
\]
\item Using the previous item we have
\begin{align*}
||\text{freq}_{k}(H_{S})-\text{freq}_{k}(G_{S})||_{1} & =\sum_{\Gamma_{S}\in\mathcal{L}_{S}(d,k)}|\text{freq}_{k}(H_{S},\Gamma_{S})-\text{freq}_{k}(G_{S},\Gamma_{S})|=\\
 & =(2t+2)\sum_{\Gamma_{S}\in\mathcal{L}_{S}(d,k)}|\sum_{\Gamma\in P_{q}(\Gamma_{S})}\text{freq}_{q}(H,\Gamma)-\sum_{\Gamma\in P_{q}(\Gamma_{S})}\text{freq}_{q}(G,\Gamma)|\le\\
 & \le(2t+2)\sum_{\Gamma_{S}\in\mathcal{L}_{S}(d,k)}\sum_{\Gamma\in P_{q}(\Gamma_{S})}|\text{freq}_{q}(H,\Gamma)-\text{freq}_{q}(G,\Gamma)|\le
\end{align*}
The double sum goes over all $\Gamma\in P_{q}(\Gamma_{S})$ with $\Gamma_{S}\in\mathcal{L}_{S}(d,k)$.
We know by \zgoto{ProjDisjointUnion} that all projection sets are
disjoint, and so each $\Gamma\in\mathcal{L}(2t+1,q)$ appears at most
once, we can therefore upper bound the double sum with the entire
set $\mathcal{L}(2t+1,q)$
\[
\le(2t+2)\sum_{\Gamma\in\mathcal{L}(2t+1,q)}|\text{freq}_{q}(H,\Gamma)-\text{freq}_{q}(G,\Gamma)|=(2t+2)\cdot||\text{freq}_{q}(H)-\text{freq}_{q}(G)||_{1}
\]
\end{enumerate}
\zSection{The Projection Subgraph}

We can look at the set $Im(T_{S})$, which contains all the simple
graphs which are transformations of $S$-graphs. An arbitrary simple
graph is not necessarily in $Im(T_{S})$, but it always contains a
subgraph which does belong to $Im(T_{S})$ (for example, the empty
subgraph). In this section, we define the projection subgraph of a
simple graph and prove that it belongs to $Im(T_{S})$. We start with
the definition.
\begin{defn}
Let $G=(V,E)$ be a simple graph. We define the \textbf{$q-$projection
subgraph} $\Psi_{q}(G)$ of $G$ as the induced subgraph of $G$ on
the vertices
\[
V\left(\Psi_{q}(G)\right)\coloneqq\bigcup_{v\in V,\text{disc}_{q}(v)\in P_{q}\left(\mathcal{L}_{S}(d,k)\right)}\text{disc}_{1}(v)
\]
\end{defn}

The idea behind the definition is that vertices in $G$ that have
$q$-discs in $P_{q}\left(\mathcal{L}_{S}(d,k)\right)$ are center-type
vertices. The $1$-disc of a center-type vertex is exactly its cluster,
and so the union of all those clusters should be in $Im(T_{S})$.
We prove this formally in the following lemma.

\zlabel{Lemma}{HardProp}
\begin{lem}
Let $G=(V,E)$ be a simple graph with maximum degree $2t+1$. Then
\end{lem}

\begin{enumerate}
\item The union that defines the vertex set $V(\Psi_{q}(G))$ is disjoint.
Namely
\[
V(\Psi_{q}(G))=\dot{\bigcup_{v\in V,\text{disc}_{q}(v)\in P_{q}\left(\mathcal{L}_{S}(d,k)\right)}}\text{disc}_{1}(v)
\]
\item The projection subgraph is in the image of $T_{S}$
\[
\Psi_{q}(G)\in Im(T_{S})
\]
\\
\end{enumerate}
\textbf{Proof} Let $v,w\in V$ be two different vertices with $\text{disc}_{q}(v),\text{disc}_{q}(w)\in P_{q}\left(\mathcal{L}_{S}(d,k)\right)$.
By definition, there are $S$-graphs $G_{S},H_{S}$ (possibly $G_{S}=H_{S}$)
and vertices $x\in V(G_{S}),y\in V(H_{S})$ such that $\text{disc}_{q}(v)=\text{disc}_{q}(x_{c})$
and $\text{disc}_{q}(w)=\text{disc}_{q}(y_{c})$. The vertices $x_{c},y_{c}$
are center-type in $G_{S},H_{S}$ with degree $2t+1$, and therefore
\[
\text{deg}_{G}(v)=\text{deg}_{G}(w)=2t+1
\]

Suppose that $\text{disc}_{1}(v)\cap\text{disc}_{1}(w)\ne\emptyset$
and in particular $w\in\text{disc}_{2}(v)$. By $\text{disc}_{q}(v)=\text{disc}_{q}(x_{c})$,
we deduce that $\text{disc}_{2}(x_{c})$ contains a vertex $z\ne x_{c}$
with $\deg_{T_{S}(G_{S})}z=2t+1$ (as ${2<3k+1=q}$). By \zgoto{TsProp}
$(3)$, we know that $z$ is center-type. Moreover, by \zgoto{DistThree}
we know that the distance between two center-type vertices in $T_{S}(G_{S})$
is at least $3$, which means $\text{dist}_{T_{S}(G_{S})}(x_{c},z)\ge3$,
in contradiction to $z\in\text{disc}_{2}(x_{c})$.

In other words, the $1$-discs of $v,w$ are disjoint, and each such
$1$-disc corresponds to the $1$-disc of a center-type vertex in
a transformation graph. These clusters can only be interpreted in
a single way, as the $2$-disc of a vertex uniquely defines its type
(\zgoto{TsProp} $(4)$). We can therefore write
\[
V(\text{disc}_{1}(v))=\{v_{in}^{1},...,v_{in}^{t},v_{out}^{1},...,v_{out}^{t},v_{c},v_{m}\}=\{v^{1},...,v^{2t},v_{c},v_{m}\}
\]
\[
V(\text{disc}_{1}(w))=\{w_{in}^{1},...,w_{in}^{t},w_{out}^{1},...,w_{out}^{t},w_{c},w_{m}\}=\{w^{1},...,w^{2t},w_{c},w_{m}\}
\]

Next, we claim that an edge between two clusters must be of the form
$(v_{in}^{i},w_{out}^{i})$ or $(v_{out}^{i},w_{in}^{i})$ for some
$i$. The vertices $v_{c},v_{m},w_{c},w_{m}$ cannot be part of inter-cluster
edges, as their counterparts in $G_{S},H_{S}$ have no inter-cluster
edges. Now, suppose $e=(v^{i},w^{j})$ is an edge. In that case, $\text{disc}_{1}(w)\subseteq\text{disc}_{4}(v)$.
However, $\text{disc}_{4}(v)=\text{disc}_{4}(x_{c})$ (as $q\ge4$)
and therefore the edge $e$ corresponds to some $e_{1}=(x^{i},z^{j})$
in $G_{S}$ between two clusters. In $G_{S}$, $\text{disc}_{4}(x_{c})$
contains the clusters of $x,z$ and so $e_{1}$ must be an edge of
the form $(x_{in}^{i},z_{out}^{j})$ or $(x_{out}^{i},z_{in}^{j})$
for some $i$. We conclude that the edge $e$ must also be between
valid end points as $\text{disc}_{4}(v)=\text{disc}_{4}(x_{c})$.

Overall, we have shown that $V\left(\Psi_{q}(G)\right)$ is a disjoint
union of clusters of size $2t+2$ that have the required ordering,
and each inter-cluster edge goes from an in-vertex to an out-vertex
with the same index. Therefore $\Psi_{q}(G)$ is exactly the transformation
graph of some $S$-graph, which means that $\Psi_{q}(G)\in Im(T_{S})$,
as required.\qed\\

Our next step is to examine the connection between the frequency vectors
of $G$ and $\Psi_{q}(G)$. Informally, if $G$ is ``close'' to
being in the image of $T_{S}$, then by \zgoto{ProjProp} $(2)$ we
would expect 
\[
\sum_{\Gamma\in P_{q}\left(\mathcal{L}_{S}(d,k)\right)}\text{freq}_{q}(G,\Gamma)\approx\frac{1}{2t+2}
\]

In that case the graph $\Psi_{q}(G)$ would be very close to being
all of $G$, and in particular the frequency distribution difference
of the two graphs should be small. We formalize this idea in the following
lemma.

\zlabel{Lemma}{ProjSubDiff}
\begin{lem}
Let $G=(V,E)$ be a simple graph with maximum degree $2t+1$. Then

\[
||\text{freq}_{q}(G)-\text{freq}_{q}(\Psi_{q}(G))||\le\left(1+2(2t+1)^{q}\right)\left(\frac{1}{(2t+2)\sum_{\Gamma\in P_{q}\left(\mathcal{L}_{S}(d,k)\right)}\text{freq}_{q}(G,\Gamma)}-1\right)
\]
\end{lem}

Note that if $\sum_{\Gamma\in P_{q}\left(\mathcal{L}_{S}(d,k)\right)}\text{freq}_{q}(G,\Gamma)\approx\frac{1}{2t+2}$
then the upper bound is close to $0$.\\

\textbf{Proof }Using \zgoto{HardProp} we have
\begin{align*}
|\Psi_{q}(G)| & =|\dot{\bigcup}_{v\in V,\text{disc}_{q}(v)\in P_{q}\left(\mathcal{L}_{S}(d,k)\right)}\text{disc}_{1}(v)|=(2t+2)|\{v\in V|\text{disc}_{q}(v)\in P\left(\mathcal{L_{S}}(d,k)\right)\}|=\\
 & =(2t+2)\sum_{\Gamma\in P_{q}\left(\mathcal{L}_{S}(d,k)\right)}\text{cnt}_{q}(G,\Gamma)=\left((2t+2)\sum_{\Gamma\in P_{q}\left(\mathcal{L}_{S}(d,k)\right)}\text{freq}_{q}(G,\Gamma)\right)|G|
\end{align*}

We can now use \zgoto{FreqDiff} to get the required frequency difference
between $G$ and its subgraph $\Psi_{q}(G)$:

\begin{align*}
||\text{freq}_{q}(G)-\text{freq}_{q}(\Psi_{q}(G))|| & \le\frac{\left(1+2(2t+1)^{q}\right)\left(|G|-|\Psi_{q}(G)|\right)}{|\Psi_{q}(G)|}=\left(1+2(2t+1)^{q}\right)\left(\frac{|G|}{|\Psi_{q}(G)|}-1\right)\\
 & =\left(1+2(2t+1)^{q}\right)\left(\frac{1}{(2t+2)\sum_{\Gamma\in P_{q}\left(\mathcal{L}_{S}(d,k)\right)}\text{freq}_{q}(G,\Gamma)}-1\right)
\end{align*}

This completes the proof of \zgoto{ProjSubDiff}.\qed

\zSection{Reduction From S-Graphs to Simple Graphs}

In this section we prove \zgoto{HardOne}. We need to show that given
an arbitrarily large $S$-graph $G_{S}$ there is a fixed size $S$-graph
$H_{S}$ with ``similar'' local structure. Our strategy is to use
the assumption of the lemma to find a fixed size approximation $H$
of $T_{S}(G_{S})$ and then use the projection subgraph $\Psi_{q}(H)$
to construct a suitable $H_{S}$.\\

\textbf{Proof of \zgoto{HardOne}} Let $d\ge2,k\ge1,\epsilon\in(0,1)$.
Let $S$ be an information set and let $A\subseteq\Omega(S)$ be a
natural set. By the assumption of the Lemma, we know that there is
a function $f$ which satisfies the required condition of \zgoto{Main}.
We define the function $f_{S,A}$ by $f_{S,A}(d,k,\epsilon)\coloneqq f(d_{1},k_{1},\epsilon_{1})$
where $d_{1}=2t+1$, $k_{1}=q$ and 
\[
\epsilon_{1}=\frac{\epsilon}{4(2t+2)^{2}\left(1+2(2t+1)^{q}\right)}
\]

Let $G_{S}=(V_{S},I)\in A$ be an $S$-graph. We need to prove that
there exists an $S$-graph $H_{S}\in A$ such that

\begin{equation}
||\text{freq}_{k}(G_{S})-\text{freq}_{k}(H_{S})||_{1}\le\epsilon\quad\wedge\quad|V(H_{S})|\le f_{S,A}(d,k,\epsilon)\label{eq:4}
\end{equation}

We start by looking at $G\coloneqq T_{S}(G_{S})$. By \zgoto{TsProp},
this is a well defined simple graph with maximal degree $d_{1}=2t+1$.
By the definition of $f$, we know that there exists a small simple
graph $H$ with
\begin{equation}
||\text{freq}_{q}(G)-\text{freq}_{q}(H)||_{1}\le\epsilon_{1}\quad\wedge\quad|V(H)|\le f(d_{1},k_{1},\epsilon_{1})\label{eq:5}
\end{equation}

Using \zgoto{ProjProp} $(2)$ for $G$, we have the following lower
bound for $H$:
\begin{align*}
\sum_{\Gamma\in P_{q}\left(\mathcal{L}_{S}(d,k)\right)}\text{freq}_{q}(H,\Gamma) & =\sum_{\Gamma\in P_{q}\left(\mathcal{L}_{S}(d,k)\right)}\text{freq}_{q}(G,\Gamma)+\sum_{\Gamma\in P_{q}\left(\mathcal{L}_{S}(d,k)\right)}\left(\text{freq}_{q}(H,\Gamma)-\text{freq}_{q}(G,\Gamma)\right)\ge\\
 & \ge\frac{1}{2t+2}-\sum_{\Gamma\in P_{q}\left(\mathcal{L}_{S}(d,k)\right)}|\text{freq}_{q}(H,\Gamma)-\text{freq}_{q}(G,\Gamma)|\ge\\
 & \ge\frac{1}{2t+2}-\sum_{\Gamma\in\mathcal{L}(2t+1,q)}|\text{freq}_{q}(H,\Gamma)-\text{freq}_{q}(G,\Gamma)|=\\
 & =\frac{1}{2t+2}-||\text{freq}_{q}(G)-\text{freq}_{q}(H)||_{1}\ge\frac{1}{2t+2}-\epsilon_{1}
\end{align*}

Using the above bound, by \zgoto{ProjSubDiff} we can bound the frequency
difference between $H$ and $\Psi_{q}(H)$:
\begin{align}
||\text{freq}_{q}(H)-\text{freq}_{q}(\Psi_{q}(H))|| & \le\left(1+2(2t+1)^{q}\right)\left(\frac{1}{(2t+2)\sum_{\Gamma\in P_{q}\left(\mathcal{L}_{S}(d,k)\right)}\text{freq}_{q}(G,\Gamma)}-1\right)\le\label{eq:6}\\
\le & \left(1+2(2t+1)^{q}\right)\left(\frac{1}{1-(2t+2)\epsilon_{1}}-1\right)\nonumber 
\end{align}

Next, we know by \zgoto{HardProp} that $\Psi_{q}(H)\in Im(T_{S})$
and so there is an $S$-graph $H_{S}^{1}$ such that $T_{S}(H_{S}^{1})=\Psi_{q}(H)$.
Finally, we know that $A$ is natural, and therefore there exists
an $S$-graph $H_{S}\in A$ with

\label{1}
\begin{equation}
||\text{freq}_{k}(G_{S})-\text{freq}_{k}(H_{S})||_{1}\le||\text{freq}_{k}(G_{S})-\text{freq}_{k}(H_{S}^{1})||_{1}\wedge|V(H_{S})|\le|V(H_{S}^{1})|\label{eq:7}
\end{equation}

We claim that this $H_{S}$ satisfies \ref{eq:4}. The size requirement
follows from \ref{eq:5}, \ref{eq:7} and \zgoto{TsProp} $(2)$

\begin{align*}
|V(H_{S})| & \le|V(H_{S}^{1})|\le|V(T_{S}(H_{S}^{1}))|=|V(\Psi_{q}(H))|\le|V(H)|\le f(d_{1},k_{1},\epsilon_{1})
\end{align*}

The frequency difference follows from \ref{eq:6}, \ref{eq:7} and
\zgoto{ProjProp} $(6)$
\begin{align*}
||\text{freq}_{k}(G_{S})-\text{freq}_{k}(H_{S})||_{1} & \le||\text{freq}_{k}(G_{S})-\text{freq}_{k}(H_{S}^{1})||_{1}\le(2t+2)\cdot||\text{freq}_{q}(T_{S}(G_{S}))-\text{freq}_{q}(T_{S}(H_{S}^{1}))||_{1}=\\
 & =(2t+2)\cdot||\text{freq}_{q}(G)-\text{freq}_{q}(\Psi_{q}(H))||_{1}\le\\
 & \le(2t+2)\cdot\left(||\text{freq}_{q}(G)-\text{freq}_{q}(H)||_{1}+||\text{freq}_{q}(H)-\text{freq}_{q}(\Psi_{q}(H))||_{1}\right)\le\\
 & \le(2t+2)\cdot\left(\epsilon_{1}+\left(1+2(2t+1)^{q}\right)\left(\frac{1}{1-(2t+2)\epsilon_{1}}-1\right)\right)\le\epsilon
\end{align*}

The last inequality is technical and does not depend on $G$. The
proof to that inequality is given in \zgoto{App1}.

Overall, we have shown that for any arbitrary $G_{S}\in A$ there
is an $S$-graph $H_{S}^{1}\in A$ of size at most ${f_{S,A}(d,k,\epsilon)=f(d_{1},k_{1},\epsilon_{1})}$,
with a local structure that approximates the local structure of $G_{S}$.
In other words, $f_{S,A}$ satisfies the condition on \zgoto{Ext}.
This completes the proof of \zgoto{HardOne}. \qed

\zSection{Proof of Winkler's Theorem for Simple Graphs}

In this section we will use the transformation $T_{S}$ and the tools
that we have developed in the previous sections to prove the simple
variant of Winkler's theorem (\zgoto{SimpleWinkler}). We will prove
this theorem by constructing a reduction from the directed edge-colored
variant. Together with the reduction from PCP to the directed edge-colored
variant (\zgoto{Winkler}), this creates a reduction from PCP to the
simple variant of Winkler's question. If we compare this to Alon's
question, then the interreducibility theorem is the ``second'' step
and what remains is to show how PCP can be reduced to the directed
edge-colored variant. Before we can prove the theorem, we need a generalized
version of the projection set, where we also consider $q$-discs of
non center-type vertices.
\begin{defn}
Let $d\ge2,k\ge1$. We define the \textbf{generalized $q-$projection
set} by
\end{defn}

\[
\widetilde{P_{q}}\coloneqq\bigcup_{G_{S}=(V_{S},I)\in A}\left\{ \text{disc}_{q}(v)|v\in V(T_{S}(G_{S}))\right\} 
\]

The set $\widetilde{P_{q}}$ contains all the possible $q$-discs
of vertices in $T_{S}(G_{S})$, including those of non center-type
vertices. In particular $P_{q}\left(\mathcal{L}_{S}(d,k)\right)\subseteq\widetilde{P_{q}}$.
We claim that if the $q$-disc set of a simple graph $G$ is contained
in the generalized $q$-projection set, then the $q$-projection subgraph
$\Psi_{q}(G)$ is the entire graph.

\zlabel{Lemma}{GenSubgraph}
\begin{lem}
Let $G=(V,E)$ be a simple graph with $\{\text{disc}_{q}(v)|v\in V(G)\}\subseteq\widetilde{P_{q}}$.
Then $G=\Psi_{q}(G)$.
\end{lem}

\textbf{Proof }The projection subgraph $\Psi_{q}(G)$ is by definition
a subgraph of $G$, so we only need to prove that any vertex $v\in G$
also belongs to $\Psi_{q}(G)$. To this end, let $v\in G$. By the
assumption on $G$, we know that $\text{disc}_{q}(v)\in\widetilde{P_{q}}$.
This means that there is a graph $H_{S}$, and a vertex $w\in V(T_{S}(H_{S}))$
with 
\begin{equation}
\text{disc}_{q}(v)=\text{disc}_{q}(w)\label{eq:8}
\end{equation}
By the definition of $T_{S}$, there is a vertex $w_{S}\in V(H_{S})$
such that $w\in cluster(w_{S})$. Thus, $w$ is either a center-type
vertex, or connected by an edge to a center-type vertex. In both cases,
$w$ belongs to the $1$-disc of some vertex with degree $2t+1$ (by
\zgoto{TsProp} $(3)$). Using \ref{eq:8}, we deduce that in $G$,
the vertex $v$ is also in the $1$-disc of a vertex with degree $2t+1$.
Denote that vertex by $z$ (possibly $z=v$). Once again, by \zgoto{TsProp}
$(3)$ we know that the $q$-disc of $z$ must correspond to a center-type
vertex in some $S$-graph and therefore $z\in\Psi_{q}(G)$. In particular
$v\in\text{disc}_{1}(z)\in\Psi_{q}(G)$, which is what we had to prove.\qed\\

The following is an immediate corollary of \zgoto{GenSubgraph} and
\zgoto{HardProp} $(2)$.

\zlabel{Corollary}{GraphInImage}
\begin{cor}
Let $G=(V,E)$ be a simple graph with $\{\text{disc}_{q}(v)|v\in V(G)\}\subseteq\widetilde{P_{q}}$.
Then $G\in Im(T_{S})$.
\end{cor}

We now have all the perquisites required to prove the main theorem
of this section.\\

\textbf{Proof of \zgoto{SimpleWinkler}} We prove the theorem by reduction
from the problem for directed edge colored graphs (i.e. $S$-graphs)
to the simple variant. To this end, let $\Phi_{S}$ be a set of $S$-graphs,
and assume that the simple problem is decidable, say by an algorithm
$Alg=Alg(d,k,\Phi)$. We will construct a deterministic algorithm
$AlgDC$ which uses $Alg$ as a subroutine to decide the directed
edge-colored variant of Winkler's problem. The main idea of the proof
is that asking the directed edge-colored question for $\Phi_{S}$
is similar to asking the simple question for the set
\[
\Phi\coloneqq\bigcup_{\Gamma_{S}\in\Phi_{S}}P_{q}(\Gamma_{S})\subseteq P_{q}\left(\mathcal{L}_{S}(d,k)\right)
\]
If $G_{S}$ is a $S$-graph whose $k$-disc set is $\Phi_{S}$, then
$T_{S}(G_{S})$ will also have $q$-discs of non center-type vertices.
To overcome this problem, we use the generalized projection set $\widetilde{P_{q}}$.
We now give the implementation of $AlgDC$.

\textbf{\uline{Algorithm}}
\begin{enumerate}
\item \textbf{function} AlgDC($d,k,\Phi_{S}=\{\Gamma_{S}^{1},...,\Gamma_{S}^{n}\}$)
\item \qquad Construct the sets $\widetilde{P_{q}}$ and $\forall i\ P_{q}(\Gamma_{S}^{i})$.
\item \qquad \textbf{for} every $\emptyset\ne X_{i}\subseteq P_{q}(\Gamma_{S}^{i})$
and $Y\subseteq\widetilde{P_{q}}\backslash P_{q}\left(\mathcal{L}_{S}(d,k)\right)$\textbf{
do}
\item \qquad\qquad \textbf{if} $Alg(t,q,\bigcup_{i=1}^{n}X_{i}\cup Y)$
is ``True'' \textbf{then}
\item \qquad\qquad \qquad\textbf{ return} ``True''
\item \qquad\qquad \textbf{end if}
\item \qquad \textbf{end for}
\item \qquad \textbf{return} ``False''
\item \textbf{end function}
\end{enumerate}
We start by explaining why the algorithm is deterministic. To construct
the sets $P_{q}(\Gamma_{S}^{i})$ and $\widetilde{P_{q}}$, we theoretically
need to go over all $S$-graphs, which cannot be done in finite time.
However, by \zgoto{SDiscSize} we know that any $q$-disc in $P_{q}(\Gamma_{S}^{i})$
can have at most $2d^{q}$ vertices, and in particular it can have
vertices of at most $2d^{q}$ clusters in the transformation graph.
This means that it is enough to examine transformation graphs with
at most $2d^{q}$ clusters, or equivalently, it is enough to examine
$S$-graphs of size at most $\frac{1}{2t+2}2d^{q}$ (by \zgoto{TsProp}
$(2)$). The for-loop runs over all possible choices of suitable $X_{i}$
and $Y$. There is only a finite amount of such choices, as $P_{q}(\Gamma_{S}^{i})$
and $\widetilde{P_{q}}$ are finite. In each iteration, the algorithm
calls $Alg$, which is deterministic by assumption, and therefore
the entire loop runs in finite time. We conclude that the entire algorithm
is deterministic, as required.\\
It remains to prove that the algorithm returns ``True'' if and only
if there is an $S$-graph $G_{S}$ with
\[
\{\text{disc}_{k}(v)|v\in V(G_{S})\}=\varPhi_{S}
\]

For the first direction, assume that there is a graph $G_{S}=(V_{S},I)$
with the $k$-disc set $\Phi_{S}=\{\Gamma_{S}^{1},...,\Gamma_{S}^{n}\}$.
If we look at the graph $G\coloneqq T_{S}(G_{S})$, then by definition
\begin{align*}
 & \{\text{disc}_{q}(v)|v\in V(G)\}=\left(\bigcup_{i=1}^{n}\left\{ \text{disc}_{q}(v_{c})|v\in V_{S},\:\text{disc}_{k}(v)=\Gamma_{S}^{i}\right\} \right)\cup\{\text{disc}_{q}(v)|v\in V(G),v\notin V_{c}\}\subseteq\\
 & \subseteq\bigcup_{i=1}^{n}P_{q}(\Gamma_{S}^{i})\cup\left(\widetilde{P_{q}}\backslash P_{q}\left(\mathcal{L}_{S}(d,k)\right)\right)
\end{align*}

Where the last inclusion is based on the fact that non center-type
vertices cannot have a $q$-disc in $P_{q}\left(\mathcal{L}_{S}(d,k)\right)$
(as even their $2$-discs differ by \zgoto{TsProp} $(4)$) . We know
by the assumption on $G_{S}$ that for any $i$ the set $X_{i}=\left\{ \text{disc}_{q}(v_{c})|v\in V_{S},\:\text{disc}_{k}(v)=\Gamma_{S}^{i}\right\} $
is not empty and so this union contains non empty subsets of $P_{q}(\Gamma_{S}^{i})$.
All the remaining $q$-discs in $G$ form a subset $Y$ of $\widetilde{P_{q}}\backslash P_{q}\left(\mathcal{L}_{S}(d,k)\right)$.
We have thus found sets $\emptyset\ne X_{i}\subseteq P_{q}(\Gamma_{S}^{i})$
and $Y\subseteq\widetilde{P_{q}}\backslash P_{q}\left(\mathcal{L}_{S}(d,k)\right)$
such that $\{\text{disc}_{q}(v)|v\in V(G)\}=\bigcup_{i=1}^{n}X_{i}\cup Y$.
In particular, the for-loop iteration that considers this choice of
$X_{i}$ and $Y$ will return ``True''.\\

For the second direction, assume that the algorithm returned ``True''.
Then, there are sets $\{X_{i}\}_{i=1}^{n}$ and $Y$ with $\emptyset\ne X_{i}\subseteq P_{q}(\Gamma_{S}^{i})$
and $Y\subseteq\widetilde{P_{q}}\backslash P_{q}\left(\mathcal{L}_{S}(d,k)\right)$
for which $Alg(t,q,\bigcup_{i=1}^{n}X_{i}\cup Y)$ returned ``True''.
By the definition of $Alg$, this means that there is a simple graph
$G$ with $\{\text{disc}_{q}(v)|v\in V(G)\}=\bigcup_{i=1}^{n}X_{i}\cup Y$
and in particular
\[
\{\text{disc}_{q}(v)|v\in V(G)\}=\bigcup_{i=1}^{n}X_{i}\cup Y\subseteq\left(\bigcup_{i=1}^{n}P_{q}(\Gamma_{S}^{i})\right)\cup\left(\widetilde{P_{q}}\backslash P_{q}\left(\mathcal{L}_{S}(d,k)\right)\right)\subseteq\widetilde{P_{q}}
\]

Note that this is the required condition of \zgoto{GraphInImage},
from which we deduce that $G\in Im(T_{S})$. In other words, there
is an $S$-graph $G_{S}$ such that $T_{S}(G_{S})=G$. We claim that
the $k$-disc set of $G_{S}$ is exactly $\varPhi_{S}$, namely
\[
\{\text{disc}_{k}(v)|v\in V(G_{S})\}=\varPhi_{S}
\]

We prove this equality by double inclusion.
\begin{enumerate}
\item Let $v\in V(G_{S})$, we claim that $\text{disc}_{k}(v)\in\varPhi_{S}$.
\\
We look at the $q$-disc $\text{disc}_{q}(v_{c})$ of the center-type
vertex corresponding to $v$ in $G$. By the definition of the $q$-projection
set, we know that 
\[
\text{disc}_{q}(v_{c})\in P_{q}(\text{disc}_{k}(v))\subseteq P_{q}\left(\mathcal{L}_{S}(d,k)\right)
\]
We also know that $\text{disc}_{q}(v_{c})\in\{\text{disc}_{q}(v)|v\in V(G)\}=\bigcup_{i=1}^{n}X_{i}\cup Y$,
but $Y\subseteq\widetilde{P_{q}}\backslash P_{q}\left(\mathcal{L}_{S}(d,k)\right)$
and so there is some index $i$ such that 
\[
\text{disc}_{q}(v_{c})\in X_{i}\subseteq P_{q}(\Gamma_{S}^{i})
\]
We conclude that $\text{disc}_{q}(v_{c})\in P_{q}(\text{disc}_{k}(v))\cap P_{q}(\Gamma_{S}^{i})$,
and so by \zgoto{ProjDisjointUnion} we get $\text{disc}_{k}(v)=\Gamma_{S}^{i}$.
By definition, $\Gamma_{S}^{i}\in\varPhi_{S}$ and so $\text{disc}_{k}(v)\in\varPhi_{S}$.
\item Let $1\le i\le n$. We claim that there is a vertex $v\in V(G_{S})$
such that $\text{disc}_{k}(v)=\Gamma_{S}^{i}$.\\
By the definition of $G$, there is a center-type vertex $v_{c}\in G$
with $\text{disc}_{q}(v_{c})\in X_{i}\subseteq P_{q}(\Gamma_{S}^{i})$.
If we look at the vertex $v\in V(G_{S})$ to which $v_{c}$ corresponds,
then by definition $\text{disc}_{q}(v_{c})\in P_{q}(\text{disc}_{k}(v))$.
We therefore have $\text{disc}_{q}(v_{c})\in P_{q}(\text{disc}_{k}(v))\cap P_{q}(\Gamma_{S}^{i})$
and so again by \zgoto{ProjDisjointUnion} we conclude that $\text{disc}_{k}(v)=\Gamma_{S}^{i}$.\\
\end{enumerate}
To conclude, we have shown that the deterministic algorithm $AlgDC(d,k,\varPhi_{S})$,
which uses $Alg$ as a subroutine, returns ``True'' if and only
if there is a $S$-graph $G_{S}$ with $\{\text{disc}_{k}(v)|v\in V(G_{S})\}=\varPhi_{S}$.
This is, of course, a contradiction to \zgoto{Winkler}. Therefore
the assumption is wrong, and there is no deterministic algorithm that
solves the simple variant of Winkler's problem, as required.\qed

\noindent \zEndChapter

\noindent \zChapter{Local Structure of Paths}

\noindent In the previous chapters we have examined the question of
finding a small graph that approximates the local structure of an
arbitrary large graph. The interreducibility theorem shows that any
model which is a natural extension of the simple model does not change
the difficulty of that question. We have conjectured that there is
a reduction from PCP to our main question (\zgoto{MainConj}). One
can think of a solution to a PCP system as a single long string, constructed
by concatenating the PCP tiles. This long string can be thought of
as an $S$-graph of a path, where each edge contains a letter or a
set of letters which comprise the string. It is therefore interesting
to ask what happens if we restrict our problem to approximating arbitrary
long paths by paths of bounded size. In this chapter, we will look
at different ways to define the ``local structure'' of a path, and
prove that in all cases the problem of finding a small approximating
path is decidable.

\zSection{Undirected Paths}

Before we examine the general case, we note that the problem of approximating
an undirected path with a small undirected path is trivial. This is
because the local structure of a path on $n$ vertices is uniquely
defined.
\begin{fact}
Let $k\ge1,\epsilon\in(0,1)$, and let $P$ be an undirected path.
Then there is an undirected path $Q$ such that
\[
||\text{freq}_{k}(P)-\text{freq}_{k}(Q)||_{1}\le\epsilon\quad\wedge\quad|V(Q)|\le\lfloor\frac{4k}{\epsilon}\rfloor+1
\]
\end{fact}

\textbf{Proof} An undirected path has a very simple local structure.
There are exactly $k$ vertices on each side of the path with $k$-discs
that include the ``boundary'' of the path, while all the rest have
the same $k$-disc $\Gamma$ which is a path of length $2k$. If we
assume that $|P|,|Q|>2k$ then
\begin{align*}
||\text{freq}_{k}(P)-\text{freq}_{k}(Q)||_{1} & =\sum_{\Delta\in\mathcal{L}(d,k)}|\text{freq}_{k}(P,\Delta)-\text{freq}_{k}(Q,\Delta)|=2k|\frac{1}{|P|}-\frac{1}{|Q|}|+|\text{freq}_{k}(P,\Gamma)-\text{freq}_{k}(Q,\Gamma)|=\\
 & =2k|\frac{1}{|P|}-\frac{1}{|Q|}|+|\frac{|P|-2k}{|P|}-\frac{|Q|-2k}{|Q|}|=4k|\frac{1}{|Q|}-\frac{1}{|P|}|
\end{align*}

If $|P|\le\lfloor\frac{4k}{\epsilon}\rfloor+1$ then taking $Q=P$
is sufficient. Otherwise we take $Q$ to be the path of length $\lfloor\frac{4k}{\epsilon}+1\rfloor$
and then

\[
||\text{freq}_{k}(P)-\text{freq}_{k}(Q)||_{1}=4k|\frac{1}{|Q|}-\frac{1}{|P|}|=4k\left(\frac{1}{|Q|}-\frac{1}{|P|}\right)\le\frac{4k}{|Q|}=\frac{4k}{\lfloor\frac{4k}{\epsilon}\rfloor+1}\le\frac{4k}{\frac{4k}{\epsilon}}=\epsilon
\]

As required.\qed

\zSection{S-Paths and S-Cycles}

The general question for directed edge-colored paths is not as trivial
as the undirected case. When we say ``directed edge-colored path''
we mean that the path graph can be written as a sequence $v_{1},...,v_{n}$
where the only edges in the path are $\forall i\ (v_{i},v_{i+1})$
and each edge attains some value $s\in S$. We start by defining this
formally.
\begin{defn}
Let $S$ be a finite information set (can represent colors, strings,
letters etc).

Let $n\in\mathbb{N}$ and let $V=\{v_{1},...,v_{n}\}$ be a set of
vertices. Let $G=(V,I)\in\Omega(S)$ be an $S$-graph.

Suppose the following holds for some $x\in\{0\}\cup\left\{ \{1\}\times S\right\} $:
\[
\forall v,w\in V\quad I(v,w)=\begin{cases}
(1,s_{i}) & \exists i\quad v=v_{i},w=v_{i+1}\\
x & v=v_{n},w=v_{1}\\
0 & \text{otherwise}
\end{cases}
\]

If $x=0$, we say that $G$ is an \textbf{$S$-path}. Otherwise, we
say that $G$ is an \textbf{$S$-cycle}. In both cases, we say that
the \textbf{size} of $G$ is $|V_{S}|=n$. 
\end{defn}

Let $k\ge1$ be an integer, and let $G=(V,I)$ be an $S$-path/$S$-cycle.
The $k$-disc of a vertex $v\in V$ is defined similarly as in the
general $S$-graph model. In addition, if $2k+2\le|G|$ then the $k$-disc
of every vertex is itself an $S$-path. Just like in the original
problem, two $k$-discs are said to be isomorphic if and only if there
is a root preserving isomorphism which also preserves $I$. We claim
that under these definitions, the problem of finding a small $S$-path
that approximates an arbitrary $S$-path is decidable.

\zlabel{Theorem}{PathByPathFormal}
\begin{thm}
Let $k\ge1,\epsilon\in(0,1)$. Let $P$ be a $S$-path, then there
is an $S$-path $Q$ such that

\[
||\text{freq}_{k}(P)-\text{freq}_{k}(Q)||_{1}\le\epsilon\ \wedge\ |Q|\le\frac{24960d^{3k}|S|^{2}L_{S}^{6}(d,k)}{\epsilon^{2}}
\]
\end{thm}

The proof is based on the following two lemmas, which will be proven
in the next sections of this chapter.

\zlabel{Lemma}{EdgeRewiring}
\begin{lem}
Let $k\ge1,\epsilon\in(0,1)$. Let $G\in\Omega(S)$ be a disjoint
union of $S$-cycles, each of size at least $2k+2$. Then, there is
an $S$-graph $H\in\Omega(S)$ which is a disjoint union of $S$-cycles
and $S$-paths such that
\end{lem}

\[
||\text{freq}_{k}(G)-\text{freq}_{k}(H)||_{1}\le\epsilon\ \wedge\ |H|\le\frac{130d^{k}|S|^{2}L_{S}^{5}(d,k)}{\epsilon}
\]

\zlabel{Lemma}{CycleBlowup}
\begin{lem}
Let $k\ge1,\epsilon\in(0,1)$. Let $G\in\Omega(S)$ be an $S$-graph
which is a disjoint union of $S$-cycles, each of size at least $2k+2$.
Then there is an $S$-path $P$ such that

\[
||\text{freq}_{k}(G)-\text{freq}_{k}(P)||_{1}\le\epsilon\ \wedge\ |P|\le\frac{8d^{k}L_{S}(d,k)}{\epsilon}|G|
\]
\end{lem}

Having stated the lemmas, we are now ready to prove the main theorem
of this chapter.\\

\textbf{Proof of \zgoto{PathByPathInformal} \textbackslash{} \zgoto{PathByPathFormal}}

Let $k\ge1,\epsilon\in(0,1)$, and let $P=(V,I)$ be an $S$-path.
We need to find an $S$-path $Q$ such that
\[
||\text{freq}_{k}(P)-\text{freq}_{k}(Q)||_{1}\le\epsilon\ \wedge\ |Q|\le\varphi
\]

Where 
\[
\varphi=\frac{24960d^{3k}|S|^{2}L_{S}^{6}(d,k)}{\epsilon^{2}}
\]

We can assume that $\varphi<|P|$ (otherwise taking $Q=P$ is sufficient),
and in particular

\[
|P|\ge\max\{2k+2,\frac{12d^{k}L_{S}(d,k)}{\epsilon}\}
\]

Denote the set of vertices of $P$ by $V(P)=\{v_{1},...,v_{n}\}$,
and let $s\in S$ be a value. We start by defining the cycle $C$
which is formed by adding the edge $I(v_{n},v_{1})=s$ to $P$. By
\zgoto{EdgeChange}, we know that
\[
||\text{freq}_{k}(P)-\text{freq}_{k}(C)||_{1}\le\frac{4d^{k}L_{S}(d,k)}{|P|}\le\frac{\epsilon}{3}
\]

We know that $2k+2<|P|=|C|$, and so in particular $C$ is a disjoint
union of $S$-cycles, each of size at least $2k+2$. By \zgoto{EdgeRewiring},
we know that there is an $S$-graph $H$ which is a disjoint union
of $S$-cycles and $S$-paths such that
\[
||\text{freq}_{k}(C)-\text{freq}_{k}(H)||_{1}\le\frac{\epsilon}{24d^{k}}\le\frac{\epsilon}{6}\ \wedge\ |H|\le\frac{3120d^{2k}|S|^{2}L_{S}^{5}(d,k)}{\epsilon}
\]

Denote the disjoint components of $H$ by $H_{1}^{C},...,H_{n_{1}}^{C}$
and $H_{1}^{P},...,H_{n_{2}}^{P}$ which are cycles and paths, respectively.
We define a new graph $\widetilde{H}$, in one of the two following
ways:
\begin{enumerate}
\item If $\sum_{i=1}^{n_{2}}|H_{i}^{P}|<2k+2$, we set $\widetilde{H}$
to be just the components $H_{1}^{C},...,H_{n_{1}}^{C}$. Each of
the vertices in $H_{1}^{P},...,H_{n_{2}}^{P}$ has a $k$-disc which
does not appear in $C$ (it is either an $S$-cycle or a small $S$-path),
and therefore
\[
\frac{\sum_{i=1}^{n_{2}}|H_{i}^{P}|}{|H|}\le\sum_{\Gamma\in\mathcal{L}_{S}(d,k)}|\text{freq}_{k}(C,\Gamma)-\text{freq}_{k}(H,\Gamma)|\le\frac{\epsilon}{24d^{k}}
\]
We then have
\begin{align*}
|H|-\sum_{i=1}^{n_{2}}|H_{i}^{P}| & \ge\frac{24d^{k}\sum_{i=1}^{n_{2}}|H_{i}^{P}|}{\epsilon}-\sum_{i=1}^{n_{2}}|H_{i}^{P}|\ge\frac{6(1+2d^{k})\sum_{i=1}^{n_{2}}|H_{i}^{P}|}{\epsilon}
\end{align*}
We use this bound, together with \zgoto{FreqDiff} to get
\begin{align*}
||\text{freq}_{k}(\widetilde{H})-\text{freq}_{k}(H)||_{1} & \le\frac{(1+2d^{k})\left(|H|-|\widetilde{H}|\right)}{|\widetilde{H}|}=\frac{(1+2d^{k})\sum_{i=1}^{n_{2}}|H_{i}^{P}|}{|H|-\sum_{i=1}^{n_{2}}|H_{i}^{P}|}\le\frac{\epsilon(1+2d^{k})\sum_{i=1}^{n_{2}}|H_{i}^{P}|}{6(1+2d^{k})\sum_{i=1}^{n_{2}}|H_{i}^{P}|}=\frac{\epsilon}{6}
\end{align*}
\item If $2k+2\le\sum_{i=1}^{n_{2}}|H_{i}^{P}|$, we set $\widetilde{H}$
to be just components $H_{1}^{C},...,H_{n_{1}}^{C}$ and also add
one addition component $H^{P}$ which is a concatenation of $H_{1}^{P},...,H_{n_{2}}^{P}$
into a single $S$-cycle where all the new edges get some arbitrary
value $s\in S$). We claim that $\widetilde{H}$ is a better approximation
of $C$ than $H$, as a result of the weight-shifting lemma (\zgoto{WeightShifting}).
Indeed, let $\Gamma\in\mathcal{L}(d,k)$ be a $k$-disc with $\text{freq}_{k}(\widetilde{H},\Gamma)<\text{freq}_{k}(H,\Gamma)$.
The only vertices in $H$ whose $k$-disc has changed are those who
were edges in $S$-paths, and so the only $k$-discs whose frequency
decreases are those which are not full $S-$paths of length $2k+1$.
There are no such vertices in $C$ (as it is a disjoint union of $S$-cycles),
and therefore $\text{freq}_{k}(C,\Gamma)=0$. By the weight-shifting
lemma, we then have
\[
||\text{freq}_{k}(\widetilde{H})-\text{freq}_{k}(C)||_{1}\le||\text{freq}_{k}(H)-\text{freq}_{k}(C)||_{1}\le\frac{\epsilon}{24d^{k}}\le\frac{\epsilon}{6}
\]
\end{enumerate}
Overall, we have constructed a graph $\widetilde{H}$ which is a disjoint
union of $S$-cycles of size at least $2k+2$. We then have by the
triangle inequality

\[
||\text{freq}_{k}(\widetilde{H})-\text{freq}_{k}(C)||_{1}\le||\text{freq}_{k}(\widetilde{H})-\text{freq}_{k}(H)||_{1}+||\text{freq}_{k}(H)-\text{freq}_{k}(C)||_{1}\le\frac{\epsilon}{6}+\frac{\epsilon}{6}\le\frac{\epsilon}{3}
\]

Finally, by \zgoto{CycleBlowup}, we know that there there is an $S$-path
$Q$ such that
\[
||\text{freq}_{k}(\widetilde{H})-\text{freq}_{k}(Q)||_{1}\le\frac{\epsilon}{3}\ \wedge\ |Q|\le\frac{8d^{k}L_{S}(d,k)}{\epsilon}|\widetilde{H}|
\]

We claim that this $S$-path $Q$ satisfies the needed requirements.
By the triangle inequality we know that the local structure of $Q$
is close to the local structure of $P$
\begin{align*}
 & ||\text{freq}_{k}(P)-\text{freq}_{k}(Q)||_{1}\le\\
 & \le||\text{freq}_{k}(P)-\text{freq}_{k}(C)||_{1}+||\text{freq}_{k}(C)-\text{freq}_{k}(\widetilde{H})||_{1}+||\text{freq}_{k}(\widetilde{H})-\text{freq}_{k}(Q)||_{1}\le\\
 & \le\frac{\epsilon}{3}+\frac{\epsilon}{3}+\frac{\epsilon}{3}=\epsilon
\end{align*}

And we also know that $Q$ is small
\[
|Q|\le\frac{8d^{k}L_{S}(d,k)}{\epsilon}|\widetilde{H}|\le\frac{8d^{k}L_{S}(d,k)}{\epsilon}\cdot\frac{3120d^{2k}|S|^{2}L_{S}^{5}(d,k)}{\epsilon}=\frac{24960d^{3k}|S|^{2}L_{S}^{6}(d,k)}{\epsilon^{2}}
\]

We can plug in $d=2$ and use the naive bound $L_{S}(d,k)=L_{S}(2,k)\le(2k)^{|S|}$
(a $k$-disc of an $S$-path can have at most $2k-1$ directional
edges that attain values in $S$) to get the bound
\[
|Q|\le24960\frac{8^{k}|S|^{2}(2k)^{6|S|}}{\epsilon^{2}}
\]

Which is an explicit bound only depending on $d,k,|S|$. This completes
the proof of the theorem.\qed\\
\\
In conclusion, we have shown that for any $\epsilon>0$ and any $S$-path
$P$, it is possible to find a small $S$-path $Q$, whose size does
not depend on $P$ such that
\[
||\text{freq}_{k}(P)-\text{freq}_{k}(Q)||_{1}\le\epsilon
\]

This means that the question of finding a small approximation is decidable
for $S$-paths, and in particular it is not possible to construct
a reduction from PCP to this variant of the problem.

\newpage
\zSection{Rewiring Edges In S-Cycles}

In this section we prove \zgoto{EdgeRewiring} by utilizing the ``rewire
and split'' technique, which was also used to prove \zgoto{FPS}.
The lemma claims that an arbitrary $S$-graph $G$ which is a disjoint
union of $S$-cycles can be approximated by a ``small'' $S$-graph
$H$ which is a disjoint union of $S$-cycles and $S$-paths. We start
by defining some additional graph related notation that will be used
in this section.\\

Let $G=(V,I)$ be an $S$-graph. Given a subset $W\subseteq V$ of
vertices, we denote by $\text{cnt}_{k}(W|G)$ the \textbf{relative
$k$-disc count vector} whose entries only count the number of $k$-discs
of each type attained by vertices in $W$. The relative $k$-disc
frequency distribution vector is defined by $\text{freq}_{k}(W|G)=\text{cnt}_{k}(W|G)/|W|$.
If $\Gamma\in\mathcal{L}_{S}(d,k)$ is a $k$-disc then $\text{cnt}_{k}(W|G,\Gamma)$
and $\text{freq}_{k}(W|G,\Gamma)$ denote the entries in the relative
vectors which correspond to $\Gamma$.

Suppose that $V_{1}\dot{\cup}V_{2}=V$ is a partitioning of $V$ into
two disjoint subsets. We define the set of all directional edges from
$V_{1}$ to $V_{2}$ by $e_{G}(V_{1},V_{2})\coloneqq|\{(x,y)|x\in V_{1},y\in V_{2},\ I(x,y)\ne0\}|$. 

The \textbf{cut} of the partition is defined by 
\[
\text{cut}_{G}(V_{1},V_{2})=e_{G}(V_{1},V_{2})+e_{G}(V_{2},V_{1})
\]
Next, we define the measure $\alpha$ of the partition:

\[
\alpha\coloneqq\alpha_{G}(V_{1},V_{2})\coloneqq\max_{\Gamma\in\mathcal{L}_{S}(d,k)}|\text{freq}_{k}(V_{1}|G,\Gamma)-\text{freq}_{k}(V_{2}|G,\Gamma)|
\]

The function $\alpha$ measures how ``balanced'' is the partition
in terms of local structure. If the local structures of $V_{1},V_{2}$
are close, then $\alpha$ would be close to $0$. Finally, we define
the value $e_{s}(P_{1},P_{2}|X,Y)$.
\begin{defn}
Let $G=(V,I)$ be a disjoint union of $S$-cycles, and let $X,Y\subseteq V$
be two subsets of $V$ (not necessarily disjoint). Given some value
$s\in S$ and two $S$-paths $P_{1},P_{2}$ of size $k$, we denote
by $e_{s}(P_{1},P_{2}|X,Y)$ the amount of $S$-paths $P=\{p_{1},...,p_{2k}\}$
in $G$ of size $2k$ where:
\end{defn}

\begin{enumerate}
\item The subgraphs induced by $\{p_{1},...,p_{k}\}$ and $\{p_{k+1},...,p_{2k}\}$
are isomorphic to $P_{1}$ and $P_{2}$, respectively.
\item The vertices $p_{k}$ and $p_{k+1}$ are in $X$ and $Y$, respectively,
with $I(p_{k},p_{k+1})=s$
\end{enumerate}
In other words, $e_{s}(P_{1},P_{2}|X,Y)$ counts the amount of edges
in $G$ that connect the two $S$-paths ($P_{1}$ and $P_{2}$) by
an edge from $X$ to $Y$ with the value $s$.\\

Our first lemma gives a basic connection between the measure $\alpha$
and the value $e_{s}(P_{1},P_{2}|X,Y)$.

\zlabel{Lemma}{MeasureConnection}
\begin{lem}
Let $k\ge1$. Let $G=(V,I)$ be a disjoint union of $S$-cycles, each
of size at least $2k+2$, and let $V_{1}\dot{\cup}V_{2}=V$ be a partitioning
of $V$. Let $s\in S$ and let $P_{1},P_{2}$ be $S$-paths of size
$k$, then
\[
|\frac{e_{s}(P_{1},P_{2}|V_{1},V)}{|V_{1}|}-\frac{e_{s}(P_{1},P_{2}|V_{2},V)}{|V_{2}|}|,|\frac{e_{s}(P_{1},P_{2}|V,V_{1})}{|V_{1}|}-\frac{e_{s}(P_{1},P_{2}|V,V_{2})}{|V_{2}|}|\le\alpha_{G}(V_{1},V_{2})|S|
\]
\end{lem}

\textbf{Proof} Suppose that $e=(x,y)$ is an edge from $x$ to $y$
that is counted by $e_{s}(P_{1},P_{2}|V_{1},V)$ or $e_{s}(P_{1},P_{2}|V_{2},V)$.
Since $x$ belongs to a cycle of size at least $2k+2$, we know that
the $k$-disc of $x$ is an $S$-path of size $2k+1$ of the form
$X(t_{1})=t_{1}P_{1}sP_{2}$ for some $t_{1}\in S$. Summing over
all possible $t_{1}$ we get
\[
e_{s}(P_{1},P_{2}|V_{1},V)=\sum_{t_{1}\in S}\text{cnt}_{k}(V_{1}|G,X(t_{1}))\quad e_{s}(P_{1},P_{2}|V_{2},V)=\sum_{t_{1}\in S}\text{cnt}_{k}(V_{2}|G,X(t_{1}))
\]

We then have
\begin{align*}
|\frac{e_{s}(P_{1},P_{2}|V_{1},V)}{|V_{1}|}-\frac{e_{s}(P_{1},P_{2}|V_{2},V)}{|V_{2}|}| & =|\frac{\sum_{t_{1}\in S}\text{cnt}_{k}(V_{1}|G,X(t_{1}))}{|V_{1}|}-\frac{\sum_{t_{1}\in S}\text{cnt}_{k}(V_{2}|G,X(t_{1}))}{|V_{2}|}|=\\
 & =|\sum_{t_{1}\in S}\text{freq}_{k}(V_{1}|G,X(t_{1}))-\sum_{t_{1}\in S}\text{freq}_{k}(V_{2}|G,X(t_{1}))\\
 & \le\sum_{t_{1}\in S}|\text{freq}_{k}(V_{1}|G,X(t_{1}))-\text{freq}_{k}(V_{2}|G,X(t_{1}))|\le\alpha_{G}(V_{1},V_{2})|S|
\end{align*}

Similarly, if $e=(x,y)$ is an edge from $x$ to $y$ that is counted
by $e_{s}(P_{1},P_{2}|V,V_{1})$ or $e_{s}(P_{1},P_{2}|V,V_{2})$,
then the $k$-disc of $y$ is an $S$-path of size $2k+1$ of the
form $P_{1}sP_{2}t_{1}$ for some $t_{1}\in S$. Similar analysis
gives us
\[
|\frac{e_{s}(P_{1},P_{2}|V,V_{1})}{|V_{1}|}-\frac{e_{s}(P_{1},P_{2}|V,V_{2})}{|V_{2}|}|\le\alpha_{G}(V_{1},V_{2})|S|
\]

This completes the proof of \zgoto{MeasureConnection}.\qed\\

We proceed by using the above lemma to give a better connection between
the measures.

\zlabel{Lemma}{MeasureConnection2}
\begin{lem}
Let $k\ge1$. Let $G=(V,I)$ be a disjoint union of $S$-cycles, each
of size at least $2k+2$, and let $V_{1}\dot{\cup}V_{2}=V$ be a partitioning
of $V$. Let $s\in S$ and let $P_{1},P_{2}$ be $S$-paths of size
$k$, then
\[
|e_{s}(P_{1},P_{2}|V_{1},V_{2})-e_{s}(P_{1},P_{2}|V_{2},V_{1})|\le2\frac{|V_{1}||V_{2}|}{|V|}\alpha_{G}(V_{1},V_{2})|S|
\]
\end{lem}

\textbf{Proof} Since $V_{1}\dot{\cup}V_{2}=V$ is a disjoint partition,
we have the following identities for $i\in\{1,2\}$
\[
e_{s}(P_{1},P_{2}|V_{i},V)=e_{s}(P_{1},P_{2}|V_{i},V_{1})+e_{s}(P_{1},P_{2}|V_{i},V_{2})
\]
\[
e_{s}(P_{1},P_{2}|V,V_{i})=e_{s}(P_{1},P_{2}|V_{1},V_{i})+e_{s}(P_{1},P_{2}|V_{2},V_{i})
\]

We use these identities to get an upper bound on the required difference

\begin{align*}
 & |e_{s}(P_{1},P_{2}|V_{1},V_{2})-e_{s}(P_{1},P_{2}|V_{2},V_{1})|=\\
 & =\frac{|V_{1}||V_{2}|}{|V|}\cdot\frac{|V_{1}|+|V_{2}|}{|V_{1}||V_{2}|}|e_{s}(P_{1},P_{2}|V_{1},V_{2})-e_{s}(P_{1},P_{2}|V_{2},V_{1})|\\
 & =\frac{|V_{1}||V_{2}|}{|V|}\cdot|\left(\frac{1}{|V_{1}|}+\frac{1}{|V_{2}|}\right)\left(e_{s}(P_{1},P_{2}|V_{1},V_{2})-e_{s}(P_{1},P_{2}|V_{2},V_{1})\right)+0-0|\\
 & =\frac{|V_{1}||V_{2}|}{|V|}\cdot|\left(\frac{1}{|V_{1}|}+\frac{1}{|V_{2}|}\right)\left(e_{s}(P_{1},P_{2}|V_{1},V_{2})-e_{s}(P_{1},P_{2}|V_{2},V_{1})\right)+\\
 & +\frac{e_{s}(P_{1},P_{2}|V_{1},V_{1})-e_{s}(P_{1},P_{2}|V_{1},V_{1})}{|V_{1}|}-\frac{e_{s}(P_{1},P_{2}|V_{2},V_{2})-e_{s}(P_{1},P_{2}|V_{2},V_{2})}{|V_{2}|}|=\\
 & =\frac{|V_{1}||V_{2}|}{|V|}\cdot|\frac{e_{s}(P_{1},P_{2}|V_{1},V_{2})}{|V_{1}|}+\frac{e_{s}(P_{1},P_{2}|V_{1},V_{2})}{|V_{2}|}-\frac{e_{s}(P_{1},P_{2}|V_{2},V_{1})}{|V_{1}|}-\frac{e_{s}(P_{1},P_{2}|V_{2},V_{1})}{|V_{2}|}+\\
 & +\frac{e_{s}(P_{1},P_{2}|V_{1},V_{1})-e_{s}(P_{1},P_{2}|V_{1},V_{1})}{|V_{1}|}-\frac{e_{s}(P_{1},P_{2}V_{2},V_{2})-e_{s}(P_{1},P_{2}|V_{2},V_{2})}{|V_{2}|}|=\\
 & =\frac{|V_{1}||V_{2}|}{|V|}\cdot|\frac{e_{s}(P_{1},P_{2}|V_{1},V_{2})+e_{s}(P_{1},P_{2}|V_{1},V_{1})}{|V_{1}|}+\frac{e_{s}(P_{1},P_{2}|V_{1},V_{2})+e_{s}(P_{1},P_{2}|V_{2},V_{2})}{|V_{2}|}\\
 & -\frac{e_{s}(P_{1},P_{2}|V_{2},V_{1})+e_{s}(P_{1},P_{2}|V_{1},V_{1})}{|V_{1}|}-\frac{e_{s}(P_{1},P_{2}|V_{2},V_{1})+e_{s}(P_{1},P_{2}|V_{2},V_{2})}{|V_{2}|}|\\
 & =\frac{|V_{1}||V_{2}|}{|V|}\cdot|\frac{e_{s}(P_{1},P_{2}|V_{1},V)}{|V_{1}|}+\frac{e_{s}(P_{1},P_{2}|V,V_{2})}{|V_{2}|}-\frac{e_{s}(P_{1},P_{2}|V,V_{1})}{|V_{1}|}-\frac{e_{s}(P_{1},P_{2}|V_{2},V)}{|V_{2}|}|\le\\
 & \le\frac{|V_{1}||V_{2}|}{|V|}\cdot\left(|\frac{e_{s}(P_{1},P_{2}|V_{1},V)}{|V_{1}|}-\frac{e_{s}(P_{1},P_{2}|V_{2},V)}{|V_{2}|}|+|\frac{e_{s}(P_{1},P_{2}|V,V_{1})}{|V_{1}|}-\frac{e_{s}(P_{1},P_{2}|V,V_{2})}{|V_{2}|}|\right)
\end{align*}

Finally, by \zgoto{MeasureConnection}, we get
\[
|e_{s}(P_{1},P_{2}|V_{1},V_{2})-e_{s}(P_{1},P_{2}|V_{2},V_{1})|\le\frac{|V_{1}||V_{2}|}{|V|}\cdot\left(\alpha_{G}(V_{1},V_{2})|S|+\alpha_{G}(V_{1},V_{2})|S|\right)=2\frac{|V_{1}||V_{2}|}{|V|}\alpha_{G}(V_{1},V_{2})|S|
\]

This completes the proof of \zgoto{MeasureConnection2}.\qed\\
\\
This result allows us to analyze our main technical tool, that is,
the rewiring of edges. We will define a condition on a partition of
$G$, and show that if that condition holds, the edges in $G$ can
be ``rewired'' in a way that decreases the cut of the partition,
without altering the FDVs of each part. We will also show that when
the condition no longer holds, the cut between the parts must be small.

\zlabel{Lemma}{RewireLemma}
\begin{lem}
Let $k\ge1,\epsilon\in(0,1)$ and let $G=(V,I_{G})$ be a disjoint
union of $S$-cycles, each of size at least $2k+2$. Suppose $V_{1}\dot{\cup}V_{2}=V$
is a partitioning of $V$. Then, either there exists an $S$-graph
$H=(V,I_{H})$ which is a disjoint union of $S$-cycles, each of size
at least $2k+2$ such that
\begin{equation}
\forall v\in V\ \text{disc}_{k}(G,v)\cong\text{disc}_{k}(H,v)\quad\wedge\quad\text{cut}_{H}(V_{1},V_{2})=\text{cut}_{G}(V_{1},V_{2})-2\label{eq:9}
\end{equation}

or the cut of $G$ is small

\begin{equation}
\text{cut}_{G}(V_{1},V_{2})\le|S|L_{S}^{2}(d,k)\left(8k+6+2\frac{|V_{1}||V_{2}|}{|V|}\alpha_{G}(V_{1},V_{2})|S|\right)\label{eq:10}
\end{equation}
\end{lem}

The proof of the lemma will be based on the following condition.

\zlabel{Condition}{RewireCond}
\begin{condition}
There exist two \textbf{disjoint isomorphic} $S$-paths of size $2k$
in $G$, denoted by $P=\{p_{1},...p_{k},p_{k+1},p_{2k}\}$ and $Q=\{q_{1},...q_{k},q_{k+1},q_{2k}\}$
such that
\begin{equation}
\text{dist}_{G}(p_{1},q_{2k}),\text{dist}_{G}(q_{1},p_{2k})\ge3\:\wedge p_{k},q_{k+1}\in V_{1}\ \wedge q_{k},p_{k+1}\in V_{2}\label{eq:11}
\end{equation}
\end{condition}

We claim that if the condition is satisfied, then Eq. (\ref{eq:9})
holds, and if not then Eq. (\ref{eq:10}) holds.\\

\textbf{Proof of \zgoto{RewireLemma}} For the first direction, suppose
there are disjoint isomorphic $S$-paths $P,Q$ of size $2k$ that
satisfy \ref{eq:11}.We define the $S$-graph $H=(V,I_{H})$ in the
following way
\[
\forall v,w\in V\quad I_{H}(v,w)=\begin{cases}
I_{G}(p_{k},p_{k+1}) & (v,w)\in\{(p_{k},q_{k+1}),(q_{k},p_{k+1})\}\\
0 & (v,w)\in\{(p_{k},p_{k+1}),(q_{k},q_{k+1})\}\\
I_{G}(v,w) & \text{otherwise}
\end{cases}
\]

In other words, we replace the edges $(p_{k},p_{k+1})$ and $(q_{k},q_{k+1})$
by $(p_{k},q_{k+1})$ and $(q_{k},p_{k+1})$. All four edges attain
the same value $s\coloneqq I_{G}(p_{k},p_{k+1})$. We observe that
\[
\text{cut}_{H}(V_{1},V_{2})=\text{cut}_{G}(V_{1},V_{2})-2
\]

as we have replaced two edges between $V_{1},V_{2}$ with two edges
which are contained within the sets. Next, we claim that $H$ is a
disjoint union of $S$-cycles of size at least $2k+2$ and that the
$k$-discs of all vertices in the graph have not changed. We distinguish
between two cases:
\begin{enumerate}
\item The $S$-paths $P,Q$ are part of the same $S$-cycle $C$ in $G$.

In this case, we can write $C=\{p_{1},...,p_{2k},x_{1},...,x_{i},q_{1},...,q_{2k},y_{1},...,y_{j}\}$.
We know by \ref{eq:11} that $\text{dist}_{G}(p_{1},q_{2k}),\text{dist}_{G}(q_{1},p_{2k})\ge3$
and therefore $i,j\ge2$. After the rewiring, $C$ becomes the two
disjoint cycles $C_{1}=\{p_{1},...,p_{k},q_{k+1},...,q_{2k},y_{1},...,y_{j}\}$
and $C_{2}=\{q_{1},...,q_{k},p_{k+1},...,p_{2k},x_{1},...,x_{i}\}$.
We then have 
\[
|C_{1}|=k+k+j\ge2k+2\quad|C_{2}|=k+k+i\ge2k+2
\]

Moreover, no $k$-discs have been affected, as the $k$-discs of all
vertices in $C_{1},C_{2}$ are the same $S$-paths of size $2k$ as
in $C$ (due to $P,Q$ being isomorphic). No other component in the
graph $G$ has been affected.
\item The $S$-paths $P$ and $Q$ are in different $S$-cycles $C_{1}$
and $C_{2}$, respectively, in $G$.

In this case, we can write $C_{1}=\{p_{1},...,p_{2k},y_{1},...,y_{j}\}$and
$C_{2}=\{q_{1},...,q_{2k},x_{1},...,x_{i}\}$ where $i,j\ge2$ (as
the $S$-cycles in $G$ are of size at least $2k+2$). After the rewiring,
the two $S$-cycles become the single $S$-cycle $C=\{p_{1},...,p_{k},q_{k+1},...,q_{2k},x_{1},...,x_{i},q_{1},...,q_{k},p_{k+1},...,p_{2k},y_{1},...,y_{j}\}$
of size $4k+i+j\ge2k+2$. Just like the first case, no $k$-discs
have changed, and no other components have been affected.
\end{enumerate}
To summarize, we have shown that if \zgoto{RewireCond} is satisfied,
then there is a graph $H$ that satisfies \ref{eq:9}.\\

In the other direction, suppose that \zgoto{RewireCond} is not satisfied.
Any $S$-path of length $2k$ can be written as a concatenation of
two $S$-paths of size $k$, connected by an edge between them. To
this end, let $s\in S$ and let $P_{1},P_{2}$ be two $S$-paths of
size $k$. We look at the $S$-path $P_{1}sP_{2}$ of size $2k$,
and examine the following sum
\[
e_{s}(P_{1},P_{2}|V_{1},V_{2})+e_{s}(P_{1},P_{2}|V_{2},V_{1})
\]

\textbf{\uline{First Case:}}\textbf{ $e_{s}(P_{1},P_{2}|V_{1},V_{2})\ne0$}

In this case, we know that there is an $S$-path $P$ in $G$ which
is isomorphic to $P_{1}sP_{2}$.

Suppose $e=(x,y)$ is an edge that has been counted by $e_{s}(P_{1},P_{2}|V_{2},V_{1})$.
By definition, the $S$-path $Q$ of size $2k$ which is formed by
taking $\text{disc}_{k}(x)\cup\text{disc}_{k}(y)$ is isomorphic to
$P$. We therefore have the $S$-paths $P,Q$ which are isomorphic
(as both are isomorphic to $P_{1}sP_{2}$). By our assumption, \zgoto{RewireCond}
is not satisfied and therefore $P,Q$ must either intersect or satisfy
\[
\text{dist}_{G}(p_{1},q_{2k}),\text{dist}_{G}(q_{1},p_{2k})<3
\]

If we look at the $S$-cycle $C=\{p_{1},...,p_{2k},c_{1},...,c_{i}\}$
to which $P$ belongs, then this requirement is equivalent to $q_{1}\in\{p_{1},...,p_{2k},c_{1},c_{2},c_{i-2k},...,c_{i}\}$,
and in particular there are at most $2k+2+(2k+1)=4k+3$ possible ways
to construct $Q$. Therefore $e_{s}(P_{1},P_{2}|V_{2},V_{1})\le4k+3$,
and so by \zgoto{MeasureConnection2}
\[
e_{s}(P_{1},P_{2}|V_{1},V_{2})+e_{s}(P_{1},P_{2}|V_{2},V_{1})\le4k+3+\left(4k+3+2\frac{|V_{1}||V_{2}|}{|V|}\alpha_{G}(V_{1},V_{2})|S|\right)
\]

\textbf{\uline{Second Case:}}\textbf{ $e_{s}(P_{1},P_{2}|V_{1},V_{2})=0$}

In this case, directly from \zgoto{MeasureConnection2} we have

\[
e_{s}(P_{1},P_{2}|V_{1},V_{2})+e_{s}(P_{1},P_{2}|V_{2},V_{1})\le0+\left(0+2\frac{|V_{1}||V_{2}|}{|V|}\alpha_{G}(V_{1},V_{2})|S|\right)
\]
We conclude that for any choice of $s\in S$ and $S$-paths $P_{1},P_{2}$
of size $k$ we have
\begin{equation}
e_{s}(P_{1},P_{2}|V_{1},V_{2})+e_{s}(P_{1},P_{2}|V_{2},V_{1})\le8k+6+2\frac{|V_{1}||V_{2}|}{|V|}\alpha_{G}(V_{1},V_{2})|S|\label{eq:12}
\end{equation}

Finally, we know that $G$ consists of $S$-cycles of length at least
$2k+2$, and so every edge in $G$ is the middle edge of exactly one
$S$-path of size $2k$. We can use this observation to define the
cut of the partition in terms of $s,P_{1},P_{2}$ in the following
way

\[
\text{cut}_{G}(V_{1},V_{2})=e_{G}(V_{1},V_{2})+e_{G}(V_{2},V_{1})=\sum_{s,P_{1},P_{2}}\left(e_{s}(P_{1},P_{2}|V_{1},V_{2})+e_{s}(P_{1},P_{2}|V_{2},V_{1})\right)
\]

Using \ref{eq:12}, and the fact that the amount of possible $k$-paths
is bounded by $L_{S}(d,k)$ (in fact even $L_{S}(d,\frac{k}{2})$)
we have
\begin{align*}
\text{cut}_{G}(V_{1},V_{2}) & =\sum_{s,P_{1},P_{2}}\left(e_{s}(P_{1},P_{2}|V_{1},V_{2})+e_{s}(P_{1},P_{2}|V_{2},V_{1})\right)=\\
 & \le|S|L_{S}^{2}(d,k)\max_{s,P_{1},P_{2}}\left(e_{s}(P_{1},P_{2}|V_{1},V_{2})+e_{s}(P_{1},P_{2}|V_{2},V_{1})\right)\le\\
 & \le|S|L_{S}^{2}(d,k)\left(8k+6+2\frac{|V_{1}||V_{2}|}{|V|}\alpha_{G}(V_{1},V_{2})|S|\right)
\end{align*}

This completes the proof of the second direction, and of \zgoto{RewireLemma}.\qed\\

We can now finally prove \zgoto{EdgeRewiring}. Given a disjoint union
of $S$-cycles $G$, our strategy would be to define a partition of
$G$ with a small $\alpha$, and then perform the edge rewiring manipulation
on the graph $G$ for as long as \zgoto{RewireCond} is satisfied.
We will then use the resulting graph to find a small approximation
of $G$.

\textbf{Proof of \zgoto{EdgeRewiring}}

Let $G$ be a disjoint union of $S$-cycles of size at least $2k+2$.
We need to find an $S$-graph $H$ which is a disjoint union of $S$-cycles
and $S$-paths such that
\[
||\text{freq}_{k}(G)-\text{freq}_{k}(H)||_{1}\le\epsilon\ \wedge\ |H|\le2\varphi
\]

Where
\[
\varphi\coloneqq\frac{65d^{k}|S|^{2}L_{S}^{5}(d,k)}{\epsilon}
\]

We can assume that $2\varphi<|G|$ (otherwise taking $H=G$ is sufficient).
We start by defining the disjoint partition $V_{1}\dot{\cup}V_{2}=V$
of $G$. For each $\Gamma\in\mathcal{L}_{S}(d,k)$, we denote the
set of vertices in $G$ with $k$-disc $\Gamma$ by $\{v_{\Gamma}^{1},...,v_{\Gamma}^{\text{cnt}_{k}(G,\Gamma)}\}$.
We then define
\[
V_{1}=\bigcup_{\Gamma\in\mathcal{L}_{S}(d,k)}\{v_{\Gamma}^{1},...,v_{\Gamma}^{\lceil\varphi\cdot\text{freq}_{k}(G,\Gamma)\rceil}\},\quad V_{2}=V/V_{1}
\]

This partition is well defined as $\lceil\varphi\cdot\text{freq}_{k}(G,\Gamma)\rceil<\lceil|G|\cdot\text{freq}_{k}(G,\Gamma)\rceil=\text{cnt}_{k}(G,\Gamma)$.
Moreover, we observe that $|V_{1}|\in\varphi\left(\varphi,\varphi+L_{S}(d,k)\right)$
as

\[
\varphi=\sum_{\Gamma\in\mathcal{L}_{S}(d,k)}\varphi\cdot\text{freq}_{k}(G,\Gamma)\le|V_{1}|\le\sum_{\Gamma\in\mathcal{L}_{S}(d,k)}\left(\varphi\cdot\text{freq}_{k}(G,\Gamma)+1\right)\le\varphi+L_{S}(d,k)
\]

We can use that to also get the bound $|V_{2}|=|V|-|V_{1}|\in(|V|-\varphi-L_{S}(d,k),|V|-\varphi)$. 

Now, let $\Gamma\in\mathcal{L}_{S}(d,k)$ be a $k$-disc. We then
have
\begin{align*}
 & |\text{freq}_{k}(V_{1}|G,\Gamma))-\text{freq}_{k}(G,\Gamma)|=|\frac{\lceil\varphi\cdot\text{freq}_{k}(G,\Gamma)\rceil}{|V_{1}|}-\text{freq}_{k}(G,\Gamma)|\le\\
 & \le\max\{\text{freq}_{k}(G,\Gamma)-\frac{\lceil\varphi\cdot\text{freq}_{k}(G,\Gamma)\rceil}{\varphi+L_{S}(d,k)},\frac{\lceil\varphi\cdot\text{freq}_{k}(G,\Gamma)\rceil}{\varphi}-\frac{\varphi\text{freq}_{k}(G,\Gamma)}{\varphi}\}\\
 & \le\max\{\frac{L_{S}(d,k)}{\varphi+L_{S}(d,k)},\frac{1}{\varphi}\}\le\frac{L_{S}^{2}(d,k)}{\varphi}
\end{align*}

And in particular, if we sum over all $\Gamma$ then
\begin{equation}
|\text{freq}_{k}(V_{1}|G))-\text{freq}_{k}(G)|\le\frac{L_{S}^{3}(d,k)}{\varphi}\label{eq:13}
\end{equation}
\\

Similar calculation can be done for $V_{2}$, if we also use the fact
that $2\varphi<|V|$:
\begin{align*}
 & |\text{freq}_{k}(G,\Gamma)-\text{freq}_{k}(V_{2}|G,\Gamma)|=\\
 & =|\text{freq}_{k}(G,\Gamma)-\frac{\text{cnt}_{k}(G,\Gamma)-\lceil\varphi\cdot\text{freq}_{k}(G,\Gamma)\rceil}{|V_{2}|}|=\\
 & =\max\{\text{freq}_{k}(G,\Gamma)-\frac{\text{cnt}_{k}(G,\Gamma)-\lceil\varphi\cdot\text{freq}_{k}(G,\Gamma)\rceil}{|V|-\varphi},\frac{\text{cnt}_{k}(G,\Gamma)-\lceil\varphi\cdot\text{freq}_{k}(G,\Gamma)\rceil}{|V|-\varphi-L_{S}(d,k)}-\text{freq}_{k}(G,\Gamma)\}\\
 & \le\max\{\frac{1}{|V|-\varphi},\frac{L_{S}(d,k)}{|V|-\varphi-L_{S}(d,k)}\}\le\max\{\frac{1}{\varphi},\frac{L_{S}(d,k)}{\varphi-L_{S}(d,k)}\}\le\frac{L_{S}^{2}(d,k)}{\varphi}
\end{align*}

By the triangle inequality, we conclude that
\[
\alpha\coloneqq\alpha_{G}(V_{1},V_{2})=\max_{\Gamma\in\mathcal{L}_{S}(d,k)}|\text{freq}_{k}(V_{1}|G,\Gamma)-\text{freq}_{k}(V_{2}|G,\Gamma)|\le\frac{2L_{S}^{2}(d,k)}{\varphi}
\]

Now, by \zgoto{RewireLemma} we know that either 
\begin{equation}
\text{cut}_{G}(V_{1},V_{2})\le|S|L_{S}^{2}(d,k)\left(8k+6+2\frac{|V_{1}||V_{2}|}{|V|}\alpha_{G}(V_{1},V_{2})|S|\right)\label{eq:14}
\end{equation}
or the graph $G$ can be replaced with a graph $\widetilde{G}$ on
the same vertex set that preserves the $k$-discs of all vertices
on one side, and has exactly two edges less between $V_{1},V_{2}$.
One can repeat this process at most a finite amount of times (as the
amount of edges between $V_{1},V_{2}$ is finite) to finally get a
graph $\widetilde{G}$ with the same $k$-discs as in $G$ which satisfies
\ref{eq:14}. In particular, we have

\begin{equation}
\text{freq}_{k}(\widetilde{G})=\text{freq}_{k}(G)\wedge\text{freq}_{k}(V_{1}|\widetilde{G})=\text{freq}_{k}(V_{1}|G)\label{eq:15}
\end{equation}

If we plug in the bounds for $|V_{1}|,|V_{2}|$ and $\alpha_{G}$,
we get
\begin{align*}
\text{cut}_{\widetilde{G}}(V_{1},V_{2}) & \le|S|L_{S}^{2}(d,k)\left(8k+6+2\frac{|V_{1}||V_{2}|}{|V|}\alpha_{G}(V_{1},V_{2})|S|\right)\le\\
 & \le|S|L_{S}^{2}(d,k)\left(8k+6+2\left(\varphi+L_{S}(d,k)\right)\frac{2L_{S}^{2}(d,k)}{\varphi}|S|\right)\le\\
 & \le|S|L_{S}^{2}(d,k)\left(8k+6+8L_{S}^{2}(d,k)|S|\right)\le16|S|^{2}L_{S}^{4}(d,k)
\end{align*}

Finally, we define the graph $H\coloneqq\widetilde{G}[V_{1}]$ which
is formed by removing all the edges between $V_{1}$ and $V_{2}$
in $\widetilde{G}$ and then taking the subgraph induced by $V_{1}$.
As a subgraph of $\widetilde{G}$ (which we know to be a disjoint
union of $S$-cycles), we know that $H$ is a disjoint union of $S$-cycles
and $S$-paths.

By \zgoto{EdgeChange}, we have a bound on the difference between
the FDVs of $H$ and $V_{1}|\widetilde{G}$

\begin{align*}
||\text{freq}_{k}(V_{1}|\widetilde{G})-\text{freq}_{k}(H)||_{1} & \le\frac{4d^{k}\left(\text{cut}_{\widetilde{G}}(V_{1},V_{2})\right)L_{S}(d,k)}{|V_{1}|}\le\frac{4d^{k}\left(16|S|^{2}L_{S}^{4}(d,k)\right)L_{S}(d,k)}{|V_{1}|}\le\frac{64d^{k}|S|^{2}L_{S}^{5}(d,k)}{\varphi}
\end{align*}

Together with \ref{eq:13} and \ref{eq:15} we conclude that
\begin{align*}
 & ||\text{freq}_{k}(G)-\text{freq}_{k}(H)||_{1}\le\\
 & \le||\text{freq}_{k}(G)-\text{freq}_{k}(V_{1}|G))||_{1}+||\text{freq}_{k}(V_{1}|G))-\text{freq}_{k}(V_{1}|\widetilde{G})||_{1}+||\text{freq}_{k}(V_{1}|\widetilde{G})-\text{freq}_{k}(H)||_{1}\le\\
 & \le\frac{L_{S}^{3}(d,k)}{\varphi}+0+\frac{64d^{k}|S|^{2}L_{S}^{5}(d,k)}{\varphi}\le\frac{65d^{k}|S|^{2}L_{S}^{5}(d,k)}{\varphi}=\epsilon
\end{align*}

And we can also bound the size of $H$
\[
|V(H)|=|V_{1}|\le\varphi+L_{S}(d,k)\le2\varphi
\]

This completes the proof of \zgoto{EdgeRewiring}. \qed

\newpage
\zSection{Blowing-Up S-Cycles}

In this section we prove \zgoto{CycleBlowup}. The lemma states that
a disjoint union of big $S$-cycles can be approximated by a single
$S$-path. The idea is that if there is a small amount of cycles,
and each cycle is very big (in terms of amount of vertices), then
``removing'' one edge from each cycle and connecting all the formed
$S$-paths together creates a single $S$-path, with a very small
amount of vertices whose $k$-discs have been affected.\\
However, if the cycles are small, than each edge removal affects a
relatively big portion of the vertices in that cycle. To handle that,
we need to start by blowing up the cycles in a way that the amount
of vertices drastically increases but the frequency does not change
by much.

The first lemma of this section shows how blowing up cycles of size
at least $2k+2$ does not affect their frequency.
\begin{lem}
Let $k\ge1$ and let $C=(V,I)$ be an $S$-cycle of size at least
$2k+2$. Denote the vertices of $C$ by $V=\{v_{1},...,v_{n}\}$ (where
each consecutive pair has an edge). Denote $I(v_{n},v_{1})$ by $c$.

Let $P$ be the $S$-path formed by removing the edge $(v_{n},v_{1})$
from $C$ and let $C_{1}$ be the $S$-cycle formed by concatenating
$m>0$ copies of $P$ one to another where the last vertex of a copy
is connected by an edge to the first vertex in the next copy and the
value of that edge is $c$. Then
\[
\text{freq}_{k}(C_{1})=\text{freq}_{k}(C)\ \wedge\ |C_{1}|=m|C|
\]
\end{lem}

\textbf{Proof }By the definition of $C_{1}$, we know that $|C_{1}|=m|P|=m|C|$.

By the way we constructed $C_{1}$, we know that each vertex in $C_{1}$
has a $k$-disc of an $S$-path which is exactly the same as its origin
vertex in $C$ (as the $S$-cycles are of size at least $2k+2$).
Moreover, if $\Gamma\in\mathcal{L}_{S}(d,k)$ is the $k$-disc of
$t$ vertices in $V$, then there are exactly $mt$ vertices in $C_{1}$
with this $k$-disc. Therefore

\begin{align*}
\text{freq}_{k}(C_{1}) & =\sum_{\Gamma\in\mathcal{L}_{S}(d,k)}\text{freq}_{k}(C_{1},\Gamma)=\frac{1}{|C_{1}|}\sum_{\Gamma\in\mathcal{L}_{S}(d,k)}\text{cnt}_{k}(C_{1},\Gamma)=\frac{1}{m|C|}\sum_{\Gamma\in\mathcal{L}_{S}(d,k)}\left(m\cdot\text{cnt}_{k}(C,\Gamma)\right)=\\
 & =\sum_{\Gamma\in\mathcal{L}_{S}(d,k)}\frac{\text{cnt}_{k}(C,\Gamma)}{|C|}=\sum_{\Gamma\in\mathcal{L}_{S}(d,k)}\text{freq}_{k}(C,\Gamma)=\text{freq}_{k}(C)
\end{align*}

This completes the proof of the lemma.\qed\\

We can now prove the Cycle Blowup lemma. Our strategy would be to
blow up the original cycles, remove a single edge from each, and finally
connect all the resulting paths into a single path.\\

\textbf{Proof of \zgoto{CycleBlowup} }Let $k\ge1,\epsilon\in(0,1)$.
Let $G\in\Omega(S)$ be an $S$-graph which is a disjoint union of
$S$-cycles, each of size at least $2k+2$. Let $t$ be the amount
of $S$-cycles in $G$. Denote this set of cycles by $C_{1},...C_{t}$.
Then
\[
\text{freq}_{k}(G)=\frac{\text{cnt}_{k}(G)}{|G|}=\sum_{i=1}^{t}\frac{\text{cnt}_{k}(C_{i})}{|G|}=\sum_{i=1}^{t}\frac{|C_{i}|}{|G|}\text{freq}_{k}(C_{i})
\]

Let 
\[
m=\lceil\frac{4d^{k}(2t-1)L_{S}(d,k)}{\epsilon|G|}\rceil\le\frac{8d^{k}tL_{S}(d,k)}{\epsilon|G|}
\]

By the previous lemma, we know that each of these cycles can be ``blown
up'' to a a cycle of size $m|C_{i}|$ such that the FDV is not changed.
Denote those new large cycles by $M_{1},...,M_{t}$. We then have
\[
\forall i\quad\text{freq}_{k}(M_{i})=\text{freq}_{k}(C_{i})\ \wedge\ |M_{i}|=m|C_{i}|
\]

Now, we remove a single edge from each of those cycles, and denote
the resulting paths by $P_{1},...,P_{t}$.

By \zgoto{EdgeChange}, we know that the removal only slightly altered
the FDVs:
\[
\forall i\quad||\text{freq}_{k}(P_{i})-\text{freq}_{k}(M_{i})||_{1}\le\frac{4d^{k}L_{S}(d,k)}{|M_{i}|}=\frac{4d^{k}L_{S}(d,k)}{m|C_{i}|}
\]

Let $H$ be the disjoint union of all the paths $P_{i}$, then
\[
\text{freq}_{k}(H)=\frac{\text{cnt}_{k}(H)}{|H|}=\sum_{i=1}^{t}\frac{\text{cnt}_{k}(P_{i})}{|H|}=\sum_{i=1}^{t}\frac{|P_{i}|}{|H|}\text{freq}_{k}(P_{i})
\]

And
\[
|H|=\sum_{i=1}^{t}|P_{i}|=\sum_{i=1}^{t}|M_{i}|=m\sum_{i=1}^{t}|C_{i}|=m|G|
\]

Finally, let $s\in S$ be some value. We define the $S$-path $P$
which is formed by connecting all the paths $P_{i}$, where the new
edges, connecting the paths, all attain the value $s$.

Once again by \zgoto{EdgeChange}, we know that $P$ is formed by
adding $t-1$ edges to $H$ and so
\begin{align*}
||\text{freq}_{k}(P)-\text{freq}_{k}(H)||_{1} & \le\frac{4d^{k}(t-1)L_{S}(d,k)}{|H|}=\frac{4d^{k}(t-1)L_{S}(d,k)}{m|G|}
\end{align*}

We claim that $P$ satisfies the required condition. First, using
the fact that $t\le|G|$ we have

\[
|P|=|H|=m|G|\le\left(\frac{8d^{k}tL_{S}(d,k)}{\epsilon|G|}\right)|G|=\frac{8d^{k}tL_{S}(d,k)}{\epsilon}\le\frac{8d^{k}L_{S}(d,k)}{\epsilon}|G|
\]

And secondly, we have
\begin{align*}
 & ||\text{freq}_{k}(P)-\text{freq}_{k}(G)||_{1}\le\\
 & \le||\text{freq}_{k}(P)-\text{freq}_{k}(H)||_{1}+||\text{freq}_{k}(H)-\text{freq}_{k}(G)||_{1}=\\
 & =||\text{freq}_{k}(P)-\text{freq}_{k}(H)||_{1}+||\sum_{i=1}^{t}\frac{|P_{i}|}{|H|}\text{freq}_{k}(P_{i})-\sum_{i=1}^{t}\frac{|C_{i}|}{|G|}\text{freq}_{k}(C_{i})||_{1}=\\
 & =||\text{freq}_{k}(P)-\text{freq}_{k}(H)||_{1}+||\sum_{i=1}^{t}\frac{m|C_{i}|}{m|G|}\text{freq}_{k}(P_{i})-\sum_{i=1}^{t}\frac{|C_{i}|}{|G|}\text{freq}_{k}(M_{i})||_{1}\le\\
 & \le||\text{freq}_{k}(P)-\text{freq}_{k}(H)||_{1}+\sum_{i=1}^{t}\frac{|C_{i}|}{|G|}|\text{freq}_{k}(P_{i})-\text{freq}_{k}(M_{i})|\le\\
 & \le\frac{4d^{k}(t-1)L_{S}(d,k)}{m|G|}+\sum_{i=1}^{t}\frac{|C_{i}|}{|G|}\frac{4d^{k}L_{S}(d,k)}{m|C_{i}|}=\\
 & =\frac{4d^{k}(t-1)L_{S}(d,k)+4d^{k}tL_{S}(d,k)}{m|G|}=\frac{4d^{k}(2t-1)L_{S}(d,k)}{m|G|}\le\epsilon
\end{align*}

We have therefore constructed a path $P$ which approximated the local
structure of $G$ and has the required size restriction. This completes
the proof of the lemma.\qed

\newpage
\zSection{Alternative Local Structure Definitions}

We have shown in the previous sections that under the standard definition
of the cnt/freq vectors, the problem of approximating an $S$-path
with a small $S$-path is decidable. Using the fact that the $k$-disc
of each vertex in a long $S$-path is an $S$-path, it is possible
to define the cnt/freq vectors differently, and then ask the question
of finding a small approximation. For example, we can count only the
left/right parts of the $k$-disc (i.e take the $S$-path that starts
at a vertex without looking ``backward'').\\
\\
In general, let $P_{k}(S)$ be the set of $k$-discs of vertices in
an $S$-path, and let $M:P_{k}(S)\to X$ be a function which maps
the $k$-discs to some finite set $X$. Then by \zgoto{FreqDiffModulo}
we have for any two $S$-paths $P,Q$
\[
||\text{freq}_{M}(P)-\text{freq}_{M}(Q)||_{1}\le||\text{freq}_{k}(P)-\text{freq}_{k}(Q)||_{1}
\]

In particular, the upper bound in \zgoto{PathByPathFormal} applies
to the corresponding function for any such $M$. \\
\\
The simplest example of such a mapping is by defining the local structure
of a vertex in an $S$-path by looking at only one ``side'' of the
$k$-disc. We give an example for the right side definition. The same
reasoning is true for the left side definition.

\zlabel{Example}{AltDef}
\begin{example}
(Right $S$-path)

Let $k\ge1$. We define the mapping function $M:P_{k}(S)\to X$ that
takes a rooted $k$-disc of a vertex in an $S$-path and returns only
the right part of that $S$-path. For example, if $P$ is an $S$-path
and $v\in V(P)$ has the $k$-disc 
\[
v_{1}\to v_{2}\to...\to v_{k}\to v\to v_{k+1}\to...\to v_{2k}
\]

then the \textbf{right $k$-disc} of $v$ is the $S$-path $v\to v_{k+1}\to...\to v_{2k}$.
We can then define the cnt/freq vectors with each entry corresponding
to possible right $k$-discs, and ask if an arbitrary $S$-path $P$
can be approximated by a small $S$-path $Q$ in terms of these new
vectors. By the above discussion, we conclude that this problem is
also decidable.\\
\end{example}

We consider a more sophisticated example. Instead of considering the
set $S$ as a set of colors, we can take it to be a finite set of
strings over the alphabet $\Sigma=\{a,b\}$. In this case, each edge
in the $S$-path represents a string, and we can think of the $k$-disc
of a vertex as a single long string that corresponds to the concatenation
of the small strings on the edges of the $k$-disc. For example, the
$3$-disc
\[
v_{1}\overset{a}{\to}v_{2}\overset{ab}{\to}v_{3}\overset{b}{\to}v\overset{aa}{\to}v_{4}\overset{bb}{\to}v_{5}\overset{aba}{\to}v_{6}
\]

corresponds to the string $aabbaabbaba$.

If we recall the Post Correspondence Problem (\zgoto{PCP}), then
a solution to a PCP system is a special pair of $S$-paths that spell
the same string. It is therefore interesting to look at the local
structure of an $S$-path where every entry corresponds to a possible
string that can be spelled by a $k$-disc.

\zlabel{Example}{AltDef2}
\begin{example}
(Concatenate edges to a string)

Let $k\ge1$ and suppose $S$ is a finite set of strings over the
alphabet $\Sigma=\{a,b\}$.

We define the mapping function $M:P_{k}(S)\to X$ that takes a rooted
$k$-disc of a vertex in an $S$-path and returns the string which
is spelled by it. For example, if $P$ is an $S$-path and $v\in V(P)$
has the $k$-disc 
\[
v_{1}\overset{s_{1}}{\to}v_{2}\overset{s_{2}}{\to}...\overset{s_{k-1}}{\to}v_{k}\overset{s_{k}}{\to}v\overset{s_{k+1}}{\to}v_{k+1}\to...\overset{s_{2k}}{\to}v_{2k}
\]

then the $k$\textbf{-string} of $v$ is the string $s_{1}s_{2}...s_{2k}$.
We can then define the cnt/freq vectors with each entry corresponding
to possible $k$-strings, and ask if an arbitrary $S$-path $P$ can
be approximated by a small $S$-path $Q$ in terms of these new vectors.
By the above discussion, we conclude that this problem is also decidable.
\end{example}

\noindent \zEndChapter

\noindent \noindent \inputencoding{latin9}\zChapter{Appendix}

\noindent \zSection{Decision Problems}

A \textbf{decision problem} is a problem that can be stated as a ``True''/''False''
question for some set of input values. A \textbf{deterministic algorithm}
is an algorithm which performs a finite amount of steps that only
depend on its input. If there is a deterministic algorithm that solves
a decision problem, then this problem is said to be \textbf{solvable}
or \textbf{decidable}. Given two decision problems $P_{1}$ and $P_{2}$,
if a deterministic algorithm that solves $P_{1}$ can also be used
as a subroutine to deterministically solve $P_{2}$, then we say that
$P_{2}$ is \textbf{reducible} to $P_{1}$.\\

In many cases, the question of decidability of a problem is reducible
to the problem of calculating a value (or a function) which depends
on the input. For example, this value can be a function representing
the amount of steps needed by an optimal Turing Machine to write the
correct output on a tape. If a well defined function can be calculated
by a deterministic algorithm, it is said to be \textbf{computable}.
If no such algorithm exists, the function is said to be \textbf{uncomputable}.\\

A very fundamental decision problem in computation theory which is
known to be undecidable was introduced by Post \zBibRef{PCP}. It
is known as the ``Post Correspondence Problem'', commonly abbreviated
as PCP.

\zlabel{Problem}{PCP}
\begin{problem}
(Post Correspondence Problem)

Let $\Sigma^{*}$ be the set of finite strings over a finite alphabet
$\Sigma$. For $s_{1},s_{2}\in\Sigma^{*}$ we denote the concatenation
of $s_{1}$ and $s_{2}$ by $s_{1}s_{2}$. A Post correspondence system
(PCS) is a \textbf{finite} set
\[
P=\left\{ (a_{1},b_{1}),...,(a_{n},b_{n})\right\} 
\]

of pairs of elements in $S$. A solution of $P$ consists of an integer
$1\le m$ and a sequence $i_{1},...,i_{m}$ such that
\[
a_{i_{1}}a_{i_{2}}...a_{i_{m}}=b_{i_{1}}b_{i_{2}}...b_{i_{m}}
\]

The Post Correspondence Problem for $P$ is to deterministically determine
whether a given set $P$ has a solution.\\
\end{problem}

The classical proof of the undecidability of PCP is by reduction from
the ``Halting Problem'' (see \zBibRef{MS} section $5.2$). It was
also shown that the problem is undecidable even when the size of $P$
is bounded \zBibRef{P5},\zBibRef{P7}.

\zlabel{Fact}{PCPUD}
\begin{fact}
PCP is undecidable.
\end{fact}

If we define the value $f(P)$ as the maximal value which needs to
be ``considered'' when searching for a solution for $P$, then this
fact is equivalent to saying that the function $f$ is not computable.
For otherwise a deterministic algorithm that computes $f(P)$ (in
a finite amount of steps) and then tries all possible sequences $i_{1},...,i_{m}$
where $m\le f(P)$ solves PCP, in contradiction to the problem being
undecidable.

\zSection{Vectors and Norms}

Suppose $v,w\in\mathbb{R}^{n}$ are two vectors, whose coordinate
representation is $v=(v_{1},...,v_{n})$ and $w=(w_{1},...,w_{n})$
accordingly. The $\ell_{1}$ norm of $v$ is defined by $||v||_{1}=\sum_{i=1}^{n}|v_{i}|$,
and the distance between $v,w$ is defined by $\text{dist}(v,w)\coloneqq||v-w||_{1}$.
For a set $W\subseteq\mathbb{R}^{n}$, the distance between $v$ and
$W$ is defined by 
\[
\text{dist}(v,W)\coloneqq||v-W||_{1}\coloneqq\inf_{w\in W}||v-w||_{1}
\]

If $W$ is finite then $\text{dist}(v,W)=\min_{w\in W}||v-w||_{1}$.
Otherwise a minimum does not necessarily exist.

\noindent \newpage
\zSection{Technical Lemmas}

\noindent \zlabel{Lemma}{App1}
\begin{lem}
Let $t,q>0$ and $\epsilon\in(0,1)$. Let $\epsilon_{1}$ be
\[
\epsilon_{1}=\frac{\epsilon}{4(2t+2)^{2}\left(1+2(2t+1)^{q}\right)}
\]
Then
\end{lem}

\[
(2t+2)\cdot\left(\epsilon_{1}+\left(1+2(2t+1)^{q}\right)\left(\frac{1}{1-(2t+2)\epsilon_{1}}-1\right)\right)\le\epsilon
\]

\noindent \textbf{Proof}

\noindent By the definition of $\epsilon_{1}$, we have 
\[
\epsilon_{1}\le\min\{\frac{1}{4t+4},\frac{\epsilon}{4t+4}\}
\]

\noindent And also
\[
2(2t+2)^{2}\cdot\left(1+2(2t+1)^{q}\right)\epsilon_{1}=\frac{\epsilon}{2}
\]

\noindent Therefore

\noindent 
\begin{align*}
 & (2t+2)\cdot\left(\epsilon_{1}+\left(1+2(2t+1)^{q}\right)\left(\frac{1}{1-(2t+2)\epsilon_{1}}-1\right)\right)=\\
 & =(2t+2)\cdot\left(\epsilon_{1}+\left(1+2(2t+1)^{q}\right)\left(\frac{(2t+2)\epsilon_{1}}{1-(2t+2)\epsilon_{1}}\right)\right)\le\\
 & \le(2t+2)\cdot\left(\frac{\epsilon}{4t+4}+\left(1+2(2t+1)^{q}\right)\left(\frac{(2t+2)\epsilon_{1}}{1-(2t+2)\frac{1}{4t+4}}\right)\right)\le\\
 & \le(2t+2)\cdot\left(\frac{\epsilon}{4t+4}+\left(1+2(2t+1)^{q}\right)\left(2(2t+2)\epsilon_{1}\right)\right)\\
 & =\frac{\epsilon}{2}+2(2t+2)^{2}\cdot\left(1+2(2t+1)^{q}\right)\epsilon_{1}=\\
 & =\frac{\epsilon}{2}+\frac{\epsilon}{2}=\epsilon
\end{align*}

\noindent Which is what we had to prove.\qed

\noindent \zEndChapter

\noindent \end{enumerate} % End the chapter enumeration

\noindent \noindent \inputencoding{latin9}\zChapterHeader{References}
\begin{itemize}[align=parleft]
\zBibSrc{AL}
{L. Lovasz. \textbf{Large Networks and Graph Limits.} American Mathematical Society, 2012.}

\noindent \zBibSrc{BER}
{P. Indyk, A. McGregor, I. Newman and K. Onak. Bertinoro workshop on sublinear algorithms 2011. \url{http://sublinear.info/42} In \textit{Open Problems in Data Streams, Property Testing, and Related Topics}, 2011.}

\noindent \zBibSrc{BSS}
{I. Benjamini, O. Schramm, and A. Shapira. Every Minor-Close Property of Sparse Graphs Is Testable. \textit{Advances in Mathematics}, 223(6):2200-2218, 2010.}

\noindent \zBibSrc{BU}
{V. K. Bulitko, \textit{On graphs with given environments of vertices.} Proc. Steklov Inst. Math. 133 (1973), 77-94}

\noindent \zBibSrc{FPS}
{H. Fichtenberger, P. Peng and C. Sohler. \textit{On Constant-Size Graphs That Preserve the Local Structure of High-Girth Graphs.} In proceedings of the 19th International Workshop on Randomization and Computation (RANDOM), 2015.}

\noindent \zBibSrc{GR}
{O. Goldreich and D. Ron. \textit{Property Testing in Bounded Degree Graphs.} Algorithmica, 32:302-343, 2002.}

\noindent \zBibSrc{HKN}
{A. Hassidim, J. Kelner, H. Nguyen, and K. Onak. Local Graph Partitions for Approximation and Testing. In \textit{Proceedings of the 50th annual IEEE sumposium on Foundations of Computer Science, pages 22-31.} IEEE, 2009}

\noindent \zBibSrc{J}
{D. P. Jacobs. \textit{Undecidability of Winkler's r-Neighborhood Problem for Covering Digraphs.} J. Combin. Theory Ser.B, 60: 254-267, 1994}\zBibSrc{MS}
{M. Sipser \textbf{Introduction to the Theory of Computation.} Cengage Learning, 2013.}

\noindent \zBibSrc{PCP}
{E. L. Post. \textit{A variant of a recursively unsolvale problem.} Bull. Amer. Math. soc. 52 (1946), 264-268}

\noindent \zBibSrc{P5}
{T. Neary \textit{Undecidability in Binary Tag Systems and the Post Correspondece Problem for Five Pair of Words.} Ins. of Neuro., Uni. of Zurich, 1998}

\noindent \zBibSrc{P7}
{Y. Matiyasevich. and G. Senizergues. \textit{Decision Problem For Semi-Thue Systems With A Few Rules.} Steklov Inst. Math. 1996}

\noindent \zBibSrc{REG}
{E. Szemeredi. Regular Partitions of Graphs. Technical report, DTIC Document, 1975.}

\noindent \zBibSrc{W}
{Peter M. Winkler. \textit{Existence of graphs with a given set of r-neighborhoods.} J. Combin. Theory Ser.B, 34(2): 165-176, 1983}

\noindent \end{itemize}

\end{document}